\documentclass[11pt]{article}
\usepackage[latin1]{inputenc}
\usepackage{amssymb}
\usepackage{amsmath}
\usepackage{amsthm}

\newcommand{\R}{\mathbb{R}}
\newcommand{\Q}{\mathbb Q}
\newcommand{\N}{\mathbb N}
\newcommand{\Z}{\mathbb Z}

\newcommand{\Fix}{\mathrm{Fix}}
\newcommand{\Per}{\mathrm{Per}}

\theoremstyle{definition}

\newtheorem{thm}{Theorem}[section]

\newtheorem{prop}[thm]{Proposition}

\newtheorem{lem}[thm]{Lemma}

\newtheorem{rem}[thm]{Remark}

\newtheorem{defi}[thm]{Definition}

\newtheorem{question}[thm]{Question}

\newcommand{\vs}{\vspace{0.3cm}}

\newcommand{\vsp}{\vspace{0.1cm}}

\date{}
\author{}

\begin{document}

\title{Rigidity for $C^1$ actions on the interval arising from hyperbolicity I:
solvable groups}
\author{C.~Bonatti, \hspace{0.1cm} I.~Monteverde, 
\hspace{0.1cm} A.~Navas \hspace{0.1cm} \& \hspace{0.1cm} C.~Rivas}
\maketitle

\vspace{-0.4cm}

\noindent{\bf Abstract:} We consider Abelian-by-cyclic groups for which the cyclic factor 
acts by hyperbolic automorphisms on the Abelian subgroup. We show that if such a group 
acts faithfully by $C^1$ diffeomorphisms of the closed interval with no global fixed point at 
the interior, then the action is topologically conjugate to that of an affine group. Moreover, 
in case of non-Abelian image, we show a rigidity result concerning the multipliers of the 
homotheties, despite the fact that the conjugacy is not necessarily smooth. Some consequences 
for non-solvable groups are proposed. In particular, we give new proofs/examples yielding the 
existence of finitely-generated, locally-indicable groups with no faithful action by $C^1$ 
diffeomorphisms of the interval.    

\vsp\vsp


\vsp\vsp


\section{Introduction}
\subsection{General panorama}

The dynamics of (non-)solvable groups of germs of diffeomorphisms around a fixed
point is an important subject that has been studied by many authors in connexion to
foliations and differential equations. There is, however, a natural group-theoretical aspect
of this study of large interest. In this direction, the classification of solvable groups of
diffeomorphisms in dimension 1 has been completed, at least in large regularity: see
\cite{burslem-wilkinson, ghys} for the real-analytic case and \cite{navas} for the 
$C^2$ case; see also \cite{bleak, navas2} for the piecewise-affine case. (For the 
higher-dimensional case, see \cite{asaoka,anne}.)

In the $C^1$ context, this issue was indirectly addresed by Cantwell and Conlon in
\cite{cantwell-conlon}. Indeed, although they were interested on problems concerning
smoothing of some codimension-1 foliations, they dealt with a particular one for which
the holonomy pseudo-group turns to be the Baumslag-Solitar group. In concrete terms, they
proved that a certain natural (non-affine) action of $BS(1,2)$ on the closed interval is non-smoothable.
Later, using the results of topological classification of general actions of $BS(1,2)$ on the interval
contained in \cite{rivas}, the whole\footnote{Some of the results of this work strongly complements
this. For instance, as we state below, the semiconjugacy is necessarily a (topological) conjugacy, 
which means that the semiconjugating map is actually a homeomorphism.}
picture was completed in \cite{guelman-liousse}: every
$C^1$ action of $BS(1,n)$ on the closed interval with no global fixed point inside is
semiconjugate to the standard affine action.

Cantwell-Conlon's proof uses exponential growth of the orbit of certain intervals to yield a
contradiction (such a behaviour is impossible close to a parabolic fixed point). This clever 
argument was later used in \cite{navas-th} to give a counter-example to the converse of the 
Thurston stability theorem: there exists a finitely-generated, locally indicable\footnote{Recall 
that a group is said to be {\em locally indicable} if every nontrivial, finitely-generated 
subgroup has a surjective homomorphism onto $\Z$. Every such group admits a 
faithful action by homeomorphisms of the interval provided it is countable; see 
\cite{fourier}.} group 
with no faithful action by $C^1$ diffeomorphisms of the interval. (See also \cite{calegari}.)
As we will see, the relation with Thurston's stability arises not only at the level 
of results. Indeed, although Cantwell-Conlon's argument is very different,
%in his thesis \cite{bonatti} (see also \cite{bonatti-annals}), the first-named author developed an
an arsenal of techniques close to Thurston's that may be applied in this context and
related ones (see {\em e.g.} \cite{anne}) was independently developed in
\cite{bonatti} (see also \cite{bonatti-annals}). The aim of this
work is to put together all these ideas (and to introduce new ones)
to get a quite complete picture of all possible
$C^1$ actions of a very large class of solvable groups, namely the Abelian-by-cyclic ones.
We will show that these actions are rigid provided the cyclic factor acts hyperbolically
on the Abelian subgroup, and that this rigidity dissapears in the non-hyperbolic case.

The idea of relating a certain notion of hyperbolicity (or at least, of growth of orbits)
to $C^1$ rigidity phenomena for group actions on 1-dimensional spaces
has been proposed --though not fully developed-- by many authors.
This is explicitly mentioned in \cite{navas-th}, while it is implicit in the examples of
\cite{shinohara}. More evidence is provided by the examples in \cite{CJN,FF, navas-growth}
relying on the original constructions of Pixton \cite{pixton} and 
Tsuboi \cite{tsuboi-pixton}. All these works suggest
that actions with orbits of (uniformly bounded) subexponential growth should be always 
$C^1$-smoothable\footnote{Actually, they should be $C^{1+\tau}$-smoothable
in case of polynomial growth, with $\tau$ depending on the degree of the polynomial;
see \cite{CJN,DKN,JNR,navas-critic}.}
(compare \cite[Conjecture 2.3]{cantwell-conlon}) and realizable in any neighborhood
of the identity/rotations \cite{navas-compositio}. Despite this evidence and
the results presented here, a complete understanding of all rigidity phenomena
arising in this context remains far from being reached. More generally,
the full picture of groups of homeomorphisms
that can/cannot act faithfully by $C^1$ diffeomorphisms remains obscure.
A particular case that is challenging from both the dynamical and the
group-theoretical viewpoints can be summarized in the next

\begin{question}
What are the subgroups of the group of piecewise affine  homeomorphisms of the
circle/interval that are topologically conjugate to groups of $C^1$ diffeomorphisms ?
\end{question}

%This will be the subject of a forthcoming work, were we plan to
%apply our methods to other classes of (non-solvable) groups.

%As a by-product, we retrieve some already-known rigidity results for $C^r$ actions of many
%groups, with $r \geq 1$, and we prove several new ones.

For simplicity, in this work, all actions are assumed to be by orientation-preserving maps.

%%%%%%%%%%%%%%%%%%%%%%%%%%%%%%%%%%%%%%%%%%%%%%%%%%%%%%%%%%%%%%%%%%%%%%%%%%%%%%%

\subsection{Statements of Results}
\label{section: statements}

Given a matrix $A=(\alpha_{i,j})\in M_d(\Z)\cap GL_d(\R)$, $d \geq 1$, 
let us consider the meta-Abelian group $G_A$ with presentation
\begin{equation} \label{eq mccarthy}
G_A := \big\langle a, b_1,\ldots, b_d\mid b_i\,b_j=b_j\, b_i\;,
\;\;ab_ia^{-1}= b_1^{\alpha_{1,i}} \cdots b_d^{\alpha_{d,i}} \big\rangle.
\end{equation}
It is known that every finitely-presented, torsion-free, Abelian-by-cyclic
group has this form \cite{BS} (see also \cite{farb mosher}).

\vs

It is quite clear that $M_d(\Z) \cap GL_d(\R) \!\subset\! GL_d(\Q)$. In particular, 
the group $G_A$ above is isomorphic to a subgroup of $\Z \ltimes_A \Q^d$.
In a slightly more general way, from now on we consider
$A\in GL_d(\Q)$ and $H$ an $A^{\pm 1}$-invariant subgroup of $\Q^d$ with 
$rank_{\Q}(H)=d$ (recall that  $rank_\Q(H)$, the $\mathbb{Q}$-rank of $H$, 
is the smallest $d'$ such that $H$ embeds into $\mathbb{Q}^{d'}$), and we let $G=\Z\ltimes_A H$.

%and we let $G \!=\! \Z \ltimes_A H$ be
%a non-Abelian finitely generated subgroup of $\Z \ltimes_A \Q^d$ such that $\Q \otimes H$ 
%has dimension (as a vector space over $\Q$) equal to $d$. This is equivalent to that 
%$rank_{\Q}(H)\! =\! d$. (Recall that $rank_{\mathbb{Q}}(H)$, the $\mathbb{Q}$-rank of $H$, 
%is the smallest $d'$ such that $H$ embeds into $\mathbb{Q}^{d'}$.) 

\vsp 

\begin{prop} \label{prop affine characterization}
{\em Suppose that the matrix $A \in GL_d (\Q)$ is $\Q$-irreducible and that
the $\Q$-rank of $H \subset \Q^d$ equals $d$. Then $\Z \ltimes_A H$ has
a faithful affine action on $\mathbb{R}$ if and only if $A$ has a positive real eigenvalue.}
\end{prop}

\vsp

Next, we assume that $A$ has all its eigenvalues of norm $\not=1$.
Our main result is the following

\vsp

\begin{thm} \label{thm main}{\em Assume $A \in GL_d(\Q)$ has no eigenvalue
of norm $1$, and let $G$ be a subgroup of $\Z \ltimes_A \Q^d$ of the
form $G = \Z \ltimes_A H$, where $rank_{\Q}(H)=d$. Then every representation
of $G$ into $\mathrm{Diff}^1_+([0,1])$ whose image group admits no global
fixed point in $(0,1)$ is topologically conjugate to a representation
into the affine group.}
\end{thm}

\vsp

For the proof of Theorem \ref{thm main}, let us begin by considering 
an action of a general group $G$ as above by homeomorphisms of 
$[0,1]$. We have the next generalization of \cite[\S 4.1]{rivas}:

\vsp

\begin{lem} \label{lem two actions} {\em Let $G$ be a group as in Theorem \ref{thm main}.
Assume that $G$ acts by homeomorphisms of the closed interval with no global fixed point
in $(0,1)$. Then either there exists $b \in H$ fixing no point in $(0,1)$, in which case the
action of $G$ is semiconjugate to that of an
affine group, or $H$ has a global fixed point in $(0,1)$, in which case the element
$a\in G$ acts without fixed points inside $(0,1)$.}
\end{lem}

\vsp

In virtue of this lemma, the proof of Theorem \ref{thm main} reduces to the next two propositions.

\vsp

\begin{prop}\label{prop no hay accione por nivels}{\em Let $G$ be a group as in Theorem
\ref{thm main}. Assume that $G$ acts by homeomorphisms of $[0,1]$
with no global fixed point in $(0,1)$. 
If the subgroup $H$ acts nontrivially but has a global fixed point inside $(0,1)$, then the action of $G$
cannot be by $C^1$ diffeomorphisms.}
\end{prop}

%\vsp

\begin{prop} \label{prop conjugacion} {\em Let $G$ be a group as in Theorem \ref{thm main}. Then every action of $G$ 
by $C^1$ diffeomorphisms of $[0,1]$ with no global fixed point in $(0,1)$
%into $\mathrm{Diff}^1_+([0,1])$ 
and having non-Abelian image is minimal on $(0,1)$.}
%In particular, any such representation is topologically conjugate to an affine action.}
\end{prop}

\vsp

The structure theorem for actions is complemented by a
result of rigidity for the multipliers of the group elements
mapping into homotheties. More precisely, we prove

\vsp

\begin{thm} \label{thm multiplier}
{\em Let $G = \Z \ltimes_A H$ be a group as in Theorem \ref{thm main},
with $a \in G$ being the generator of $\mathbb{Z}$ (whose action on $H$ is given by $A$).
Assume that $G$ acts by $C^1$ diffeomorphisms of $[0,1]$ with no fixed point in $(0,1)$ 
and the image group is non-Abelian. Then the derivative
of $a$ at the interior fixed point coincides with the ratio of the homothety to which $a$ is
mapped under the homomorphism of $G$ into the affine group given by Theorem
\ref{thm main}. More generally, for each $k \neq 0$ and all $b \in H$, the derivative
of $a^k b$ at its interior fixed point equals the $k^{th}$-power of the ratio of
that homothety.}
\end{thm}

\vsp

Besides several consequences of the preceding theorem given in the next section, there 
is an elementary one of particular interest. Namely, if we consider actions as in Theorem 
\ref{thm main} but allowing the possibility of global fixed points in $(0,1)$, then only 
finitely many components of the complement of the set of these points are such that the 
action restricted to them has non-Abelian image. Otherwise, the element $a$ would admit 
a sequence of hyperbolic fixed points, all of them with the same multiplier, converging 
to a parabolic fixed point, which is absurd. 

Another consequence of the previous results concerns centralizers. For simplicity, we just 
give an statement involving the Baumslag-Solitar's group, yet a more general version 
certainly holds for groups as in Theorem \ref{thm main}.

\vsp

\begin{prop} {\em The centralizer inside $\mathrm{Diff}_+^1([0,1])$ of a subgroup $G$ isomorphic to 
$BS(1,2) := \langle g,h \!: ghg^{-1} = h^2 \rangle$ is the group of diffeomorphisms having support in 
the complement of the support of $h$. In particular, if $G$ has no global fixed point in the interior, 
then its centralizer is trivial.}
\label{centers}
\end{prop}

\vsp

%\noindent{\em Proof:} 
Indeed, let $I$ be a closed  interval in $[0,1]$ restricted to which the action of $h$ has no global fixed 
point in the interior. We need to show that every element $f$ of the centralizer of $G$ fixes $I$ 
and acts trivially on it. To do this, first notice that, by Theorem \ref{thm main}, the group $G$ 
fixes $I$, and its action on it is topologically conjugate to that of an affine group. 
%Moreover, $f$ commutes with $h$ and all its powers (including the fractional ones). 
Therefore, if $f$ fixes $I$, then by commutativity it must fix the unique fixed point of $g$ inside $I$. 
Again by commutativity, the set of fixed points of $f$ is $G$-invariant, and since the $G$-orbits inside $I$ 
are dense, we conclude that $f$ acts trivially on $I$.
%to it must be contained in the corresponding translation flow. In particular, it must act trivially on $I$, otherwise it wouldn't commute with $g$. 
If $f$ does not fix $I$, then it moves it into a disjoint interval, so that $I, f(I), f^2 (I), \ldots$ are infinitely many pairwise disjoint 
intervals restricted to which the $G$-action has non-Abelian image, which was shown to be impossible just after the 
statement of Theorem \ref{thm multiplier}.
%By Theorem \ref{thm multiplier} and commutativity we have $Dg (h^n(a)) = Dg (a) = 2$ for all $n \geq 1$, which is impossible as 
%the sequence $h^n (a)$ should accumulate at a parabolic fixed point of $g$. This proves Proposition \ref{centers}.
%$\hfill\square$
 
\vsp\vsp\vsp

%\begin{prop} {\em Let $f$ be a $C^1$ diffeomorphism of $[0,1]$ that centralizes 
%the action of a subgroup of $\mathrm{Diff}_+^1([0,1])$ isomorphic to 
%$BS(1,2) := \langle g,h \!: ghg^{-1} = h^2 \rangle$. Then the interior 
%of the supports of $f$ and $h$ are disjoint. In particular, $f$ is the identity 
%when restricted to certain open intervals.}
%\end{prop}

%Indeed, by Theorem \ref{thm main}, on any interval $I$ where $h$ acts 
%with no interior fixed points yet fixing the endpoints, the action of $BS(1,2)$ is 
%topologically conjugate to that of an affine group. Recall that $f$ commutes with 
%$h$ and all its powers (including the fractional ones). On the one hand, if $f$ fixes 
%$I$, then its restriction to it must be contained in the corresponding translation 
%flow. In particular, it must act trivially on $I$, otherwise it wouldn't commute 
%with $g$. On the other hand, if $f$ does not fix $I$, then it moves into a disjoint 
%interval, so that $I, f(I), f^2 (I), \ldots$ are pairwise disjoint. Letting $a$ be the 
%fixed point of $g$ on $I$, by Theorem \ref{thm multiplier} and commutativity 
%we have $Dg (h^n(a)) = Dg (a) = 2$ for all $n \geq 1$, which is impossible as 
%the sequence $h^n (a)$ should accumulate at a parabolic fixed point of $g$. 

%\vsp\vsp

Theorem \ref{thm multiplier} could lead one to think that the topological conjugacy to the affine
action is actually smooth at the interior.\footnote{In general, the conjugacy above is not
smooth at the endpoints even in the real-analytic case: see \cite{burslem-wilkinson} for 
a very complete discussion on this.} (Compare \cite{shub-sullivan}.) Nevertheless, a 
standard application of the Anosov-Katok technique leads to $C^1$ (faithful) actions 
for which this is not the case. As we will see, in higher regularity, the rigidity holds: 
if $r \geq 2$, then for every faithful action by $C^r$ diffeomorphisms with no interior 
global fixed point, the conjugacy is a 
$C^r$ diffeomorphism at the interior. It seems to be an interesting problem to 
try to extend this rigidity to the class $C^{1+\tau}$. Another interesting problem 
is to construct actions by $C^1$ diffeomorphisms that are conjugate to actions of 
non-Abelian affine groups though they are non-ergodic with respect to the Lebesgue 
measure. (Compare \cite{KM}.)

\vsp

The hyperbolicity assumption for the matrix $A$ is crucial for the validity of Theorem
\ref{thm main}. Indeed, Abelian groups of diffeomorphisms acting nonfreely (as those
constructed in \cite{tsuboi-pixton}) provide easy counter-examples with all
eigenvalues equal to 1. A more delicate construction leads to the next

\vsp

\begin{thm}\label{actuan}
{\em Let $A \in GL_d (\mathbb{Q})$ be non-hyperbolic and $\Q$-irreducible. Then
$G := \Z \ltimes_A \Q^d$ admits a faithful action by $C^1$ diffeomorphisms
of the closed interval that is not semiconjugate to an affine one though has
no global fixed point in $(0,1)$.}
\end{thm}

\vsp

This work is closed by some extensions of our main theorem to actions by 
$C^1$ diffeomorphisms 
of the circle. Roughly, we rule out Denjoy-like actions in class $C^1$ for the 
groups $G$ above, though such actions may arise in the continuous cathegory 
(and also in the Lipschitz one; see \cite[Proposition 2.3.15]{navas-book}).  
In particular, we have:

\vsp

\begin{thm}\label{thm circulo}
{\em Let $G$ be a group as in Theorem \ref{thm main}. Assume that $G$ acts 
by $C^1$ diffeomorphisms of the circle with non-Abelian image. Then the 
action admits a finite orbit.}
\end{thm}

\vsp

This theorem clarifies the whole picture. Up to a finite-index subgroup $G_0$, the action has
global fixed points. The group $G_0$ can still be presented in the form $\Z \ltimes_{A^k} H_0$
for a certain $k \geq 1$; as $A^k$ is hyperbolic, and application of Theorem \ref{thm main} to
the restriction of the action of $G_0$ to intervals between global fixed points shows that
these are conjugate to affine actions. Thus, roughly, $G$ is a finite (cyclic) extension of a
subgroup of a product of affine groups acting on intervals with pairwise disjoint interior. 
Moreover, only finitely many of these affine groups can be non-Abelian. (Otherwise, by 
Theorem \ref{thm multiplier}, there would be accumulation of hyperbolic fixed points 
of $a^k$ with the same multiplier towards a parabolic fixed point.)

To conclude, let us mention that the examples provided by Theorem \ref{actuan} can be
adapted to the case of the circle. More precisely, if $A \in GL_d(\Q)$ is non-hyperbolic
and $\Q$-irreducible, then $G := \Z \ltimes_A \Q^d$ admits a faithful action by $C^1$
circle diffeomorphisms having no finite orbit.

%%%%%%%%%%%%%%%%%%%%%%%%%%%%%%%%%%%%%%%%%%%%%%%%%%%%%%%%%%%%%%%%%%%

\subsection{Some comments and complementary results/examples}
\label{sec extensions}

Although the results presented so far only concern certain solvable groups, they lead to 
relevant results for other classes of groups. Let us start with an almost direct consequence 
of Theorem~\ref{thm multiplier}. For any pair of positive integers $m,n$, let $BS(1,m;1,n)$ 
be the group defined by 
$$BS(1,m;1,n) := \big\langle a,b,c \mid aba^{-1}=b^m \!, \hspace{0.05cm} 
aca^{-1}=c^n \big\rangle = BS(1,m) *_{\langle a \rangle} BS(1,n).$$  
%(in particular, the image of the subgroup generated by $b,c$ is rank-2 free)
In other words, the subgroups generated by $a,b$  and $a,c$ are isomorphic 
to $BS(1,m)$ and $BS(1,n)$,  respectively,  and no other relation is assumed.

Notice that every element $\omega \in BS(1,m;1;n)$ can be written in a unique way as 
$\omega = a^k \omega_0$, where $k\!\in\!\Z$ and $\omega_0$ belongs to the subgroup 
generated by $b,c$ and their roots. One easily deduces that $BS(1,m;1,n)$ is locally 
indicable, hence it admits a faithful action by homeomorphisms of the interval (see the 
second footnote in page 1). However, it is easy to give a more explicit embedding of 
$BS(1,m;1,,n)$ into $\mathrm{Homeo}_+([0,1])$. Indeed, start by associating to $a$ 
a homeomorphism $f$ without fixed points in $(0,1)$. Then  choose a fundamental 
domain $I$ of $f$ and homeomorphisms $g_0,h_0$ supported on $I$ and 
generating a rank-2 free group. Finally, extend $g_0$ and $h_0$ into 
homeomorphisms $g,h$ of $[0,1]$ so that $fgf^{-1}=g^m$ and 
$fhf^{-1}=h^n$ hold. Then the action of $B(1,m;1,n)$ defined 
by associating $g$ to $b$ and $h$ to $c$ is faithful. 

In what concerns smooth actions of $BS(1,m;1,n)$ on the interval, we have: 

\vsp 

\begin{thm}\label{c.2BS} {\em Let $m,n$ be distinct positive integers. Given a 
representation of $B(1,m;1,n)$ into $\mathrm{Diff}^1_+([0,1])$, let us denote by 
$f,g,h$ the images of $a,b,c$, respectively. Then, the interiors of the supports of $g$ 
and $h$ are disjoint. In particular, $g$ and $h$ commute, hence the action is not faithful. }
\end{thm}

\noindent{\em Proof:} The supports of $g$ and $h$ consist of unions of segments bounded 
by successive non-repelling fixed points of $f$; in particular, any two of these segments either 
coincide or have disjoint interior. If one of these segments is contained in the support of 
$g$ (resp., $h$), then Theorem \ref{thm multiplier} asserts that its interior contains a 
unique hyperbolically-repelling fixed point of $f$ with derivative equal to $m$ (resp., $n$). 
Since $m\neq n$, the open segments in the supports of $g$ and $h$ must be disjoint. 
$\hfill\square$
 
\vsp
\vsp

\begin{rem} Theorem \ref{c.2BS} admits straightforward generalizations replacing the 
Baumslag-Solitar groups $BS(1,m)$ and $BS(1,n)$ by groups associated (as in Theorem 
\ref{thm main}) to matrices $A$ and $B$ that are hyperbolically expanding ({\em i.e.} 
with every eigenvalue of modulus $>1$) and have distinct eigenvalues. 
\end{rem}

\vsp

Below we give two other results in the same spirit. The first of these is 
new, whereas the second is already known though our methods
provide a new and somewhat more conceptual proof. More 
sophisticated examples will be treated elsewhere.

Let us consider the group $G_{\lambda,\lambda'}$ generated by the transformations
of the real-line
$$c \!: x\mapsto x+1,
\quad
b \!: x \mapsto \lambda x,
\quad
a \!: x \mapsto sgn(x) |x|^{\lambda'},$$
where $\lambda, \lambda'$ are positive numbers. These groups are known to be
non-solvable for certain parameters $\lambda'$. Indeed, if $\lambda'$ is a
prime number, then the elements $a$ and $c$ generate a free group (see \cite{cohen}).

\vsp

\begin{thm} {\em For all integers $m,n$ larger than 1, the group
$G_{m,n}$ does not embed into the group $C^1$
diffeomorphisms of the closed interval.}
\end{thm}

\noindent{\em Proof:} Assume that $G_{m,n}$ can be realized as a group of $C^1$ 
diffeomorphisms of $[0,1]$. Then Theorem \ref{thm main} applies to both subgroups 
$\langle b,c \rangle$ and $\langle a,b \rangle$ (which are isomorphic to $BS(1,m)$ and 
$BS(1,n)$, respectively). Let us consider a maximal open subinterval $I = ( x_0,x_1 )$ 
that is invariant under $c$ and where the dynamics of $c$ has no fixed point. The 
relation $bcb^{-1} = c^m$ shows that the action of $b$ on $I$ is nontrivial. Proposition 
\ref{prop no hay accione por nivels} then easily implies that $b$ preserves $I$, and by 
Theorem \ref{thm main}, the restriction of the action of $\langle b,c \rangle$ to $I$ 
is conjugate to an affine action. Let $y$ be the fixed point of $b$ inside $I$. As 
before, the relation $aba^{-1} = b^n$ forces $a$ to fix all points $x_0,y,x_1$; moreover, 
the actions of $\langle a,b \rangle$ on both intervals $(x_0,y)$ and $(y,x_1)$ are 
conjugate to affine actions. Finally, notice that the relation $aba^{-1} = b^n$ forces 
the derivative of $b$ to be equal to 1 at $y$. However, this contradicts Theorem 
\ref{thm multiplier} when applied to $\langle b,c \rangle$. $\hfill\square$

\vsp
\vsp
\vsp

As another application of our results, we give an alternative proof of a theorem 
from \cite{navas-th}:

\vsp

\begin{thm} {\em If $\Gamma$ is a non-solvable subgroup of $SL_2(\R)$, then
$\Gamma \ltimes \Z^2$ does not embed into $\mathrm{Diff}_+^1 ([0,1])$.}
\label{th-th}
\end{thm}

\noindent{\em Proof:} Since $\Gamma$ is non-solvable, it must contain 
two hyperbolic elements $A,B$ generating a free group. Theorem \ref{thm main}
applied to $\Z \ltimes_A \Z^2 \subset \Gamma \ltimes \Z^2$ implies that the action
restricted to $\langle A,\Z^2 \rangle$ is topologically conjugate to an affine action with
dense translation part on each connected component $I$ fixed by $\langle A,\Z^2 \rangle$ 
and containing no point that is globally fixed by this subgroup. 
As $B$ normalizes $\Z^2$, it has to be affine in the coordinates induced by this 
translation part. As a consequence, the action of $\Gamma \ltimes \Z^2$ is that 
of an affine group on each interval $I$ as above. We thus conclude that the action 
factors throughout that of a solvable group, hence it is unfaithful. $\hfill\square$

\vsp\vsp

\begin{rem}
It is not hard to extend the previous proof to show that $\Gamma \ltimes \Z^2$ does 
not embed into the group of $C^1$ diffeomorphisms of neither the open interval nor 
the circle. (Compare \cite[\S 4.2]{navas-th} and \cite[\S 4.3]{navas-th}, respectively.)
\end{rem}

\begin{rem} \label{th-counter}
All groups discussed in this section are locally indicable. We thus get 
different infinite families of finitely-generated, locally-indicable groups with no 
faithful actions by $C^1$ diffeomorphisms of the closed interval. The existence 
of such groups was first established in \cite{navas-th}; the examples considered 
therein correspond to those of Theorem \ref{th-th}. 
\end{rem}

%%%%%%%%%%%%%%%%%%%%%%%%%%%%%%%%%%%%%%%%%%%%%%%%%%%%%%%%%%%%%%%%%%%

\section{On affine actions}
\label{section-remark-base-canonica}

Before passing to proofs, let us make an important remark. So far, in all statements concerning groups of the 
form $\mathbb{Z} \rtimes_A H$, the $\Q$-rank of $H$ is assumed to be maximal. Therefore, we may fix a 
$\Q$-basis $b_1,\ldots,b_d$ of $\Q^d$ made of elements in $H$. 
%Up to multiplying these vectors by a large-enough integer, we can (and will) assume that they 
%all correspond to elements in $H$, which will be denoted in the same way. \marginpar{no entiendo}
More importantly, 
since all conditions to be imposed on $A$ (if any) are invariant under conjugacy, up to changing 
$A$ by a conjugate matrix in $GL_d (\Q)$, we may assume that $b_1,\ldots,b_d$ is the canonical 
basis of $\Q^d$. This assumption will be made in order to simplify specific computations.

In this section, we prove Proposition \ref{prop affine characterization}. To simplify, vectors
$(t_1,\ldots,t_d) \in \R^d$ will be denoted horizontally, though must be viewed as 
vertical ones. We begin with

\vsp

\begin{prop}\label{propafin}{\em Given $A \in GL_d(\Q)$, let $G$ be a subgroup
of $\Z \ltimes_A \Q^d$ of the form $\Z \ltimes_A H$, where $H$ contains the canonical basis 
$\{ b_1 \ldots, b_d\}$ of $\Q^d$ (so that, in particular, $rank_{\Q}(H) = d$).

\noindent (i) If $(t_1,...,t_d)\in\R^d$ is an eigenvector of the transpose
$A^T$ with eigenvalue $\lambda\in\R_+\!\setminus\!\{1\}$, then there exists
a homomorphism $\psi:G\rightarrow \mathrm{Aff}_+(\R)$ with
non-Abelian image defined by $\psi(b_i) := T_{t_i}$ and $\psi(a) := M_\lambda$,
where $T_t$ and $M_{\lambda}$ stand for the translation by an amplitude $t$
and the multiplication by a factor $\lambda$, respectively. This
homomorphism is injective if and only if $\{t_1,...,t_d\}$
is a $\Q$-linearly-independent subset of $\R$.\\
\noindent (ii) Every homomorphism
$\psi\!: G \rightarrow \mathrm{Aff}_+(\R)$ with non-Abelian
image is conjugate to one as those described in (i).}
\end{prop}

\noindent {\em Proof:} The first claim of item (i) is obvious. For the other claim, 
notice that injectivity of $\psi$ is equivalent to injectivity of its restriction to $H$. 
%As explained at the beginning  of \S \ref{section: statements} ????, we may fix a basis 
%Let $\{b_1,\ldots,b_d\}$ be a basis of $\Q^d$ contained in $H$. 
We let $a$ be the generator of the $\Z$-factor of $G$. Assume there is an 
element  $b \in H$  mapping into the trivial translation. This element writes as 
$b = b_1^{\beta_1} \cdots b_d^{\beta_d} \in H$ for certain $\beta,\ldots,\beta_d$ in $\Q$. 
Then we have $\sum_i \beta_i t_i = 0$, which implies that the $t_i$'s are linearly dependent over
$\Q$. Conversely, assume $\sum_i \beta_i t_i = 0$ holds for certan rational numbers $\beta_i$ that are not
all equal to zero. Up to multiplying them by a common integer factor, we may assume that such a
relation holds with all $\beta_i$'s integer. Then $b := b_1^{\beta_1} \cdots b_d^{\beta_d}$
is a nontrivial element of $H$ mapping into the trivial translation under $\psi$.

For (ii), 
%we further assume (with no loss of generality) that $b_1,\ldots,b_d$ is the canonical basis of $\Q^d$. S
suppose $\psi:G \rightarrow \mathrm{Aff}_+(\R)$ is a homomorphism
with non-Abelian image. Then we have
$$\{ id \} \subsetneq \psi([G,G]) \subseteq[\mathrm{Aff}_+(\R),\mathrm{Aff}_+(\R)]
= \{T_t,\ t \in\R\}.$$
Fix $b\in[G,G]$ such that $\psi(b)$ is a nontrivial translation.
As $b\in H$, we have that $\psi(b)$ commutes with every element in
$\psi(H)$. Therefore, $\psi(H)$ is a subgroup of the group of translations.

Let $t_1,\ldots,t_d$ in $\mathbb{R}$ be such that $\psi(b_i)=T_{t_i}$. As
$\psi(G)$ is non-Abelian, we have $\psi(a)=T_t M_\lambda$ for certain
$\lambda\neq 1$ and $t \in \mathbb{R}$. We may actually suppose that $t=0$
just by conjugating $\psi$ by $T_{\frac t {\lambda-1}}$. Then, if we choose 
$k \in \mathbb{N}$ so that all $k \alpha_{i,j}$ belong to $\Z$, for each
$i \! \in \! \{1,\ldots, d\}$ we have 
$$T_{\lambda k t_i}
= \psi(a)\psi(b_i^k)\psi(a)^{-1}
= \psi \big( (a b_i a^{-1} )^k \big)
= \psi(b_1^{k\alpha_{1,i}}\ldots b_d^{k\alpha_{d,i}})
= T_{k \alpha_{1,i}t_1+\ldots + k \alpha_{d,i}t_d}.$$
Thus, $\lambda t_i=\alpha_{1,i} t_1+\ldots+\alpha_{d,i} t_d$, which
implies that $(t_1,\ldots,t_d)$ is an eigenvector of $A^T$ with eigenvalue
$\lambda$. $\hfill\square$

\vspace{0.1cm}

\begin{rem} The preceding proposition implies in particular that if $A^T$ has no real eigenvalue,
then there is no representation of $G$ in $\mathrm{Aff}_+(\R)$ with non-Abelian image. As a
consequence, due to Theorem \ref{thm main}, if moreover the eigenvalues of $A^T$ all have
modulus different from 1, then every representation of $G$ in $\mathrm{Diff}_+^1([0,1])$
has Abelian image.

As a matter of example, given positive integers $m,n$, let $A$ be the matrix
$$A=A_{m,n}:=\left(\begin{array}{cc}
m & -n\\
n & m
\end{array}\right).$$
Then the group $G(m,n) := \Z \ltimes_A \Q^2$ has no inyective representation into
$\mathrm{Diff}^1_+([0,1])$. Notice that each of these groups $G(m,n)$ is locally indicable. 
Hence, this produces still another infinite family of finitely-generated, locally-indicable
groups with no faithful action by $C^1$ diffeomorphisms of the closed interval.  
(Compare Remark \ref{th-counter}.)
\end{rem}

In view of the discussion above, the proof of Proposition \ref{prop affine characterization}
is closed by the next

\begin{lem}
 {\em Suppose that the matrix $A \in GL_d(\Q)$ is $\Q$-irreducible. If $\lambda\in\R$
is an eigenvalue of $A^T$ and $v := (t_1,\ldots,t_d)\in\R^d$ is such that
$A^T v = \lambda v$, then $\{t_1,\ldots,t_d\}$ is a $\Q$-linearly-independent
subset of $\R$.}
\end{lem}

\noindent {\em Proof:} If $v$ is an eigenvector of $A^T$, then the
subspace $v^\perp \subseteq\R^d$ is invariant under $A$. Since $A$ is
$\Q$-irreducible, we have $v^\perp\cap \Q^d = \{  0 \}$. Therefore,
if $v := (t_1,\ldots,t_d)$ and $\beta_1,\ldots,\beta_d$ in $\Q$ 
verify $\beta_1 t_1 + \cdots + \beta_d t_d=0$, then we have 
$\beta_1=\ldots=\beta_d=0$. $\hfill\square$

%%%%%%%%%%%%%%%%%%%%%%%%%%%%%%%%%%%%%%%%%%%%%%%%%%%%%%%%%%%%

\section{On continuous actions on the interval}

%%%%%%%%%%%%%%%%%%%%%%%%%%%%%%%%%%%%%%%%%%%%%%%%%%%%%%%%%%%%

In this section, we deal with actions by homeomorphisms. The proof
below was given in \cite{rivas} for the Baumslag-Solitar group $B(1,2)$.
As we next see, the argument can be adapted to the group $G$.

\vsp
\vsp
\vsp

\noindent {\em Proof of Lemma \ref{lem two actions}:} Assume that $H$ has a global fixed point $x$.  
Then the $a$-orbit of $x$ is made up of global fixed points of $H$. Therefore, this $a$-orbit has to 
accumulate at both $0$ and $1$, which means that $a$ has no fixed point in $(0,1)$. Indeed, otherwise  
this $a$-orbit would accumulate into a common fixed point $x_* \in (0,1)$ of both $H$ and $a$, which 
would mean that $x_*$ is $G$-invariant, which is against our assumption.

We next show that if $G$ acts by homeomorphisms of $[0,1]$ 
in such a way that $H$ admits no global fixed point on $(0,1)$, then 
the action is semiconjugate to that of an affine group. To do this, we let 
$N\subseteq H$ be the set of elements having a fixed point inside $(0,1)$. 
As $H$ is Abelian, $N$ is easily seen to be a subgroup. 
We claim that $N$ is strictly contained in $H$. Indeed, let 
$\{ b_1,\ldots, b_d \} \subset H$ be a $\mathbb{Q}$-basis of $\mathbb{Q} \otimes H$. We affirm that 
one of these generators $b_i$ has no fixed point. Indeed, 
if $b_1$ has no fixed point, then we are done. Otherwise, let $x_1$ be a fixed point of $b_1$. If $b_2$ 
has no fixed point, then we are done. Otherwise, let $x_2$ be a fixed point of $b_2$. If $x_2$ 
is fixed by $b_1$, then we have found a common fixed point $x_* := x_2$. If not, then either $b_1 (x_2)$  
or $b_1^{-1} (x_2)$ is closer to $x_1$ than $x_2$. In the former case, let $x_* := \lim_{n \to \infty} b_1^n (x_2)$, 
and in the latter case, let $x_*:= \lim_{n \to -\infty} b_1^n (x_2)$. By commutativity, we have that $x_*$ is fixed by 
both $b_1$ and $b_2$ in each case. Now repeating this argument finitely many times, we will detect a common 
fixed point for all the elements $b_i$ provided each of them has a fixed point. Nevertheless, this common fixed 
point will be fixed by all rational powers of the generators, hence by all elements in $H$, which is contrary 
to our assumption.

%Since $H$ has finite $\Q$-rank, we have $N\neq H$. (This easily follows along the lines of \cite[Exercise 2.2.47]{navas-book} 
%just noticing that every homeomorphism of the interval has the same fixed points as each of its rational powers.)
We now claim that there is an $H$-invariant infinite measure $\nu$ on $(0,1)$ that is finite on compact subsets.  
(Actually, this follows from \cite[Proposition 2.2.48]{navas-book}, but we repeat the argument here because is so simple.) 
Indeed, let $h \in H$ be an element having no fixed point in $(0,1)$. Then $H / \langle h \rangle$ naturally acts on the space 
$(0,1) / \! \! \sim$, where the equivalence relation $\sim$ corresponds to being in the same orbit under $h$. 
As $h$ is fixed-point free, $(0,1) / \!\! \sim$ is topologically a circle. Since $H / \langle h \rangle$ is Abelian, 
its action on this topological circle preserves a probability measure. This measure lifts to a measure 
$\nu$ on $(0,1)$ that is finite on compact sets and is invariant under $H$.

We next claim that $\nu$ has no atoms, and that it is unique up to scalar multiple. Indeed, by \cite{plante} 
(see also \cite[\S 2.4.5]{navas}), this holds whenever  
$H/N$ is not isomorphic to $\Z$, and this is the case here because $N$ is $A^T$-invariant and $A^T$ has no eigenvalue of modulus 1. 

%Now, following \cite{plante}, this implies that $\nu$ cannot have atoms. Indeed, $H / N$ freely acts on the totally 
%ordered set $Fix (N)$, hence by H\"older's theorem

Now, as $H$ is normal in $G$, we have that $a_*(\nu)$ is also invariant
by $H$. (Remind that, by definition, $a_*(\mu) (S):=\mu (a^{-1}(S))$.) By uniqueness, 
this implies that $a_*(\nu)=\lambda\nu$ for some
$\lambda\in\R_{>0}$. More generally, for every $g \in G$, there exists
$\lambda_g \in\R_{>0}$ such that $g_*(\nu) = \lambda_g \nu$. The map
$g \mapsto \lambda_g$ is a homomorphism from $G$ into $\R^*_{>0}$.
It is then easy to check that the map $\psi \! :G\to \mathrm{Aff}_+(\R)$, $g\mapsto \psi_g$, defined by
$$\psi_g(x) := \frac{1}{\lambda_b}x + \nu \big( [1/2,b(1/2)] \big)$$
is a representation, where $\nu ([q,p]) := -\nu([p,q])$ for 
$q > p$. Indeed, for $f,g\in G$, we have
\begin{eqnarray*}
\psi_{fg}(x)
&=&\frac{1}{\lambda_f}\left(\frac{1}{\lambda_g} x + \mu \big( [f^{-1}(1/2),g(1/2)] \big) \right)\\ 
&=& \frac{1}{\lambda_f}\left(\frac{1}{\lambda_g} x +\mu \big( [1/2,g(1/2)] \big) + \mu \big( [f^{-1}(1/2),1/2] \big) \right)\\
&=& \frac{1}{\lambda_f} \left(\frac{1}{\lambda_g}x+\mu \big( [1/2,g(1/2)] \big) \right) + \mu \big( [1/2,f(1/2)] \big)\\ 
&=& \psi_f\circ \psi_g(x).\end{eqnarray*}
Moreover, the map $F \!:\R\rightarrow\R$ defined 
by $F(x):=\nu([1/2,x])$ semiconjugates the action of $G$ with $\psi$. Indeed, 
\begin{eqnarray*}
F (g(x))
&=& \mu \big( [1/2,g(x)] \big)\\
&=& \mu \big( [g(1/2),g(x)] \big) +\mu \big( [1/2,g(1/2)] )\\
&=& g^{-1}_*\mu \big( [1/2,x] \big) + \mu \big( [1/2,g(1/2)] \big)\\
&=& \frac{1}{\lambda_g} F(x) + \mu \big( [1/2,g(1/2)] \big)\\
&=& \psi_g \big( F(x) \big),
\end{eqnarray*}
as desired.
$\hfill\square$\\

\vsp

In the statement of Lemma \ref{lem two actions}, the semiconjugacy is not 
necessarily a conjugacy. This easily follows by applying a Denjoy-like technique 
replacing the orbit of a single point by that of a wandering interval. See also Theorem
\ref{teo semiconj} below, where this procedure is carried out smoothly on the open interval.

%%%%%%%%%%%%%%%%%%%%%%%%%%%%%%%%%%%%%%%%%%%%%%%%%%%%%%%%%%%%

\section{On $C^1$ actions on the interval}

\subsection{All actions are semiconjugate to affine ones}
\label{seccionprueba}

%%%%%%%%%%%%%%%%%%%%%%%%%%%%%%%%%%%%%%%%%%%%%%%%%%%%%%%%%%%%

In this section, we show Proposition \ref{prop no hay accione por nivels}. Suppose 
for a contradiction that for an action of $G$ by $C^1$ diffeomorphisms of $[0,1]$ 
without global fixed points in $(0,1)$, the subgroup $H$ acts nontrivially
on $(0,1)$ but having a fixed point inside. For each $x\in (0,1)$ which is not fixed by $H$,
let us denote by $I_x$ the maximal interval containing $x$ such that $H$ has no fixed point
inside. Since $G$ has no global fixed point in $(0,1)$ and $H$ is normal in $G$, we must
have \hspace{0.01cm} $a^n(I_x)\cap I_x =\emptyset \quad \text{for all} \quad n\not = 0.$
\hspace{0.1cm} In particular, $I_x$ is contained in $(0,1)$. Moreover, $a$ has no fixed point in
$(0,1)$, and up to changing it by its inverse, we may suppose that $a(z) > z$ for all $z \in (0,1).$ 
Notice that the set of fixed points of $H$ is invariant under $a$, hence it accumulates at both 
endpoints of $[0,1]$.

The rough idea now is, for a point $x$ not fixed by $H$ as above, to apply $a^{-1}$ iteratively 
at $x$ and examine the behavior of an appropriately defined {\em displacement vector} (see 
(\ref{displacement-vector}) below). Our first lemmata (\ref{lem tauxconA} to \ref{triconA}) build the groundwork needed 
to show that the direction of this vector nearly converges along a subsequence (Lemma \ref{tauestri}). That 
$A$ is hyperbolic then implies that, along this subsequence, the magnitude of the vector is 
uniformly expanded (Lemma \ref{lem proyectivo}), giving a contradiction.

To implement the strategy above, we first recall a useful tool that arises in this context, namely, 
there is an $H$-invariant infinite measure $\mu_x$ supported on $I_x$ which is finite
on compact subsets. This measure is not unique, but independently of the choice,
we can define the {\em translation number homomorphism} $\tau_{\mu_x} \!: H\to \R$
by $\tau_{\mu_x}(h) := \mu_x([z,h(z)))$   
%(here and in all what follows, we use the convention $\mu([y,z)) := -\mu((z,y])$ for $z < y$). 
(This value is independent of $z \in I_x$.)  The kernel $K_x$ of this homomorphism coincides with the set 
of elements of $H$ having fixed points inside $I_x$; see \cite[Section 2.2.5]{navas-book} for all of this.

From now on, as explained at the beginning of \S \ref{section-remark-base-canonica}, 
we let $\{b_1,\ldots,b_d\}$ be the canonical basis 
of $\Q^d$, which (with no loss of generality) we assume to be contained in $H$. Although unnatural, 
this choice equips $\R\otimes H$ with an inner product, which yields to the following crucial notion.

\begin{defi} For every $I_x$ as above, we define the {\em translation vector}
$\vec{\tau}_{\mu_x} \in \R^d$ as the unit vector (with respect to the max norm) 
pointing in the direction $(t_1,\ldots, t_d)$, where $t_i := \mu_x ([z, b_i(z)))$.
\end{defi}

%The homomorphism $\phi_{\preceq_x}$ can be
%recovered from the measure $\mu_x$. Indeed, the map $h \mapsto \mu_x ([y, h(y)))$
%is a group homomorphism from $H$ into $\mathbb{R}$ (the value is independent of
%$y \in I_x$), and there is a constant $\lambda>0$ such that
%\begin{equation}\label{eq unique morphism}\phi_{\preceq_x}(h)
%=\lambda \, \mu_x([x, h(x))) \text{, for all } h\in H.\end{equation}

%\begin{lem} The direction $\tau_{\mu_x}$ is perpendicular to the hyperplane
%$ker()$.
%\end{lem}

%\begin{lem}\label{lem base} {\em If we identify each $b = b_1^{\beta_1} \cdots b_d^{\beta_d} \in H$ with the
%vector $(\beta_1,\ldots,\beta_d) \in \R \otimes H$, then $\vec{\tau}_{\mu_x} \in \mathbb{R}^d$
%is orthogonal to the hyperplane $\R\otimes K_x$.}
%\end{lem}

%\noindent{\em Proof:} For $b = b_1^{\beta_1} \cdots b_d^{\beta_d} \in H$ we have that $\tau_{\mu_x} (b)
%= \sum_i \beta_i \mu_x ([z,b_i(z)))$ equals zero if and only if the vector $(\beta_1,\ldots,\beta_d)$
%is ortoghonal to $\vec{\tau}_{\mu_x}$. $\hfill\square$

%\vs
\vs

In the sequel, we will denote $\vec{\tau}_{\mu_x}$
simply by $\vec{\tau}_x$. We have

\begin{lem}\label{lem tauxconA}
%{\em For each $x \in (0,1)$, we have $K_{a^{-1}(x)}=A^{-1}K_x$. Moreover, 
{\em The directions of $\vec{\tau}_{a^{-1}(x)}$ and $A^T \vec{\tau}_x$ coincide.}
%(where $A^T$ is the transpose of $A$).
\label{lem:vary}
\end{lem}

\noindent{\em Proof:} %A vector $v=(\beta_1,...,\beta_d) \in \R \otimes H$ 
%gives a positive (resp., negative) value under $\tau_{\mu_{a^{-1}(x)}}$ if and only if 
%$$b_1^{\beta_1} \cdots b_d^{\beta_d}(a^{-1}(z)) > a^{-1}(z)$$
%(resp., $b_1^{\beta_1} \cdots b_d^{\beta_d}(a^{-1}(z)) < a^{-1}(z)$) holds for every $z \in I_x$, that is
%$$a^{-1}b_1^{\alpha_{1,1} \beta_1+ \cdots +\alpha_{1,d} \beta_d}\cdots
%b_d^{\alpha_{d,1} \beta_1+ \cdots +\alpha_{d,d} \beta_d}(z) > a^{-1}(z)$$
%(resp., $a^{-1}b_1^{\alpha_{1,1} \beta_1+ \cdots +\alpha_{1,d} \beta_d}\cdots
%b_d^{\alpha_{d,1} \beta_1+ \cdots +\alpha_{d,d} \beta_d}(z) < a^{-1}(z)$). This happens if and only if
%$$b_1^{\alpha_{1,1} \beta_1+ \cdots +\alpha_{1,d} \beta_d}\cdots b_d^{\alpha_{d,1} \beta_1+ \cdots
%+\alpha_{d,d} \beta_d}(z) > z$$
%(resp., $b_1^{\alpha_{1,1} \beta_1+ \cdots +\alpha_{1,d} \beta_d}\cdots b_d^{\alpha_{d,1} \beta_1+
%\cdots +\alpha_{d,d} \beta_d}(z)<z$).
%This directly yields the first assertion of the lemma. 
%The second one is an easy consequence. Indeed, 
%$\vec{\tau}_{a^{-1}x}$ generates the subspace
%$$(\R \otimes K_{a^{-1}x})^{\perp} =
%(A^{-1} (\R \otimes K_x))^\perp = A^T (\R \otimes K_x)^\perp,$$
%and the last subspace is generated by the vector $A^T\vec{\tau}_x$.
Since $a^*(\mu_x)$ (remind that $a^{*} (\nu) (S) := \nu (a(S))$) is an $H$-invariant measure 
on $I_{a^{-1}(x)}$, by definition, we have that the $i^{th}$ entry of $\vec{\tau}_{a^{-1}(x)}$ 
coincides with $a^{*} (\mu_x) \big( [a^{-1}(x), b_i (a^{-1}(x)) \big).$ 
Thus, this entry equals
\begin{eqnarray*}
\mu_x \big( [x, ab_ia^{-1}(x)] \big)
&=& 
\tau_{\mu_x} (a b_i a^{-1})\\
&=&
\frac{1}{k} \tau_{\mu_x} (a b_i^k a^{-1})\\
&=& 
\frac{1}{k} \tau_{\mu_x} (b_1^{k \alpha_{1,i}} b_2^{k \alpha_{2,i}} \cdots b_d^{k \alpha_{i,d}}).
\end{eqnarray*}
Choosing $k$ so that all $k \alpha_{i,j}$ belong to $\mathbb{Z}$, this yields
$$\mu_x \big( [x, ab_ia^{-1}(x)] \big) 
= \frac{1}{k} \sum_{j = 1}^{d} k \alpha_{j,i} \tau_{\mu_x} (b_j) 
= \sum_{j=1}^d \alpha_{j,i} \vec{\tau}_x (b_j),$$
as desired.
$\hfill\square$

\vs

We now state our main tool to deal with $C^1$ diffeomorphisms. Roughly, it says that
diffeomorphisms that are close-enough to the identity in the $C^1$ topology 
behave like translations under
composition. For each $\delta > 0$, we denote $U_{\delta}(id)$ the neighborhood of
the identity in $\mathrm{Diff}_+^1([0,1])$ given by
$$U_{\delta}(id) := \Big\{ f \in \mathrm{Diff}_+^1([0,1])
\!: \sup_{z \in [0,1]} \big| Df(z) - 1 \big| < \delta \Big\}.$$

\begin{prop}[\cite{bonatti}]
\label{prop bonatti} {\em For each $\eta>0$ and all $k\in\N$, there exists a neighborhood $U$ of the identity
in $\mathrm{Diff}^1_+([0,1])$ such that for all $f_1,\ldots , f_k$ in $U$, all $\epsilon_1,\ldots , \epsilon_k$
in $\{-1,1\}$ and all $x\in [0,1]$, we have
$$\Big| [f_k^{\epsilon_k}\circ \ldots \circ f_1^{\epsilon_1}(x) - x]
- \sum_i \epsilon_i(f_i(x)-x) \Big| \leq \eta \max_j \big\{ |f_j(x)-x| \big\}.$$}
\end{prop}

\noindent{\em Proof:}
First of all, observe that if $g\in \mathrm{Diff}^1_+([0,1])$ satisfies $|Dg(z)-1|<\lambda$
for all $z\in[0,1]$, then for all $x,y$,
$$\big| (g(x)-x)-(g(y)-y) \big| < \lambda|x-y|.$$

Next, notice that for every $f \in U_{\delta}(id)$ and all $x \in [0,1]$, 
there exists $y \in [0,1]$ such that 
\begin{small}
$$\big| (f_i^{-1}(x) - x) - (x - f_i(x)) \big|
=
\big| (f_i(x) - x) - (f_i(f_i^{-1}(x)) - f_i^{-1}(x)) \big|\\
=
\big| Df_i(y) - 1 \big| \cdot \big| x - f_i^{-1}(x) \big|\\ 
\leq
\delta |x - f_i^{-1}(x)|.$$
\end{small}Using this, it is not hard to see that we may assume that $\epsilon_i=1$ for all $i$. 

We proceed by induction on $k$.
The case $k=1$ is trivial. Suppose the lemma holds up to $k-1$, and choose $\delta>0$ so that
the lemma applies to any $k-1$ diffeomorphisms in the neighborhood $U=U_{\delta}(id)$ for the
constant $\eta/2$. We may suppose $\delta$ is small enough to verify $\delta(k-1 +\eta/2)<\eta/2$.
Now take $f_1,\ldots, f_k$ in $U_{\delta}(id)$ and $x\in[0,1]$. Then the value of the expression
$$\Big| f_k\circ \ldots\circ f_1 (x)-x -\displaystyle\sum_{i=1}^k (f_i(x)-x) \Big|$$
is smaller than or equal to
$$\Big| f_k\circ \ldots\circ f_1 (x)- f_{k-1}\circ \ldots\circ f_1 (x) -(f_k(x)-x) \Big|
+ \Big| f_{k-1}\circ \ldots\circ f_1 (x)-x -
\displaystyle\sum_{i=1}^{k-1} (f_i(x)-x) \Big|.$$
Now notice that, by the inductive hypothesis, the second term in the sum above
is bounded from above by \hspace{0.03cm} $\eta/2 \max_j|f_j(x)-x|$. 
\hspace{0.03cm} Moreover, the observation at the beginning of the
proof and the inductive hypothesis yield
$$\big| f_k(f_{k-1}\circ \ldots\circ f_1 (x))- f_{k-1}\circ \ldots\circ f_1 (x)-(f_k(x)-x) \big|
\leq
\delta \hspace{0.02cm} \big| f_{k-1}\circ \ldots\circ f_1 (x)-x \big|\\
\leq
\delta \Big( \sum_{i=1}^{k-1}|f_i(x)-x|+\varepsilon \Big),$$
with $\varepsilon< \eta/2 \max_j|f_j(x)-x|$. By the choice of $\delta$,
the last expression is bounded from above by  $\eta/2 \max_j|f_j(x)-x|$,
thus finishing the proof. $\hfill\square$

\vs
\vs

%\vs

%\begin{lem} \label{lem powers}
%{\em Given $\eta > 0$ and $q \in \mathbb{N}$, there exists $\delta > 0$ 
%satisfying the following: if $f \in \mathrm{Diff}_+^1([0,1])$ is such that
%$f^{1/q} \in U_{\delta}(id)$, then for all $x \in [0,1]$,
%$$\big| (f(x)-x) - q(f^{1/q}(x)-x)\big| \leq \eta \hspace{0.04cm} |f^{1/q}(x)-x|.$$
%In particular,}
%$$\frac{| f(x) - x |}{q + \eta} \leq |f^{1/q} (x) - x| \leq \frac{|f(x) - x|}{q - \eta}.$$
%\end{lem}

%\noindent{\em Proof:} This directly follows from the proposition by letting \hspace{0.01cm} 
%$k = q$ \hspace{0.01cm} and \hspace{0.01cm} $f_1 = \ldots = f_k = f^{1/q}$. 
%\hspace{0.01cm} Details are left to the reader. 
%$\hfill\square$

%\vs
%\vs

%Notice that the lemma above does not state that if $f$ is close to the identity, then its roots
%(whenever they exist) remain close to the identity. Indeed, this is known to be false in general.
%Nevertheless, as we will see along the proof of the lemma below, this turns to be partially true in the group $G$. 
We next deduce some consequences from this proposition. 
To do this, first recall that the set of global fixed points of $H$ accumulate 
at both endpoints of $[0,1]$. Hence, given an element $b \in H$, for each $\delta > 0$, there is 
$\sigma_1 > 0$ which is  fixed by $H$ and such that $b$ restricted to $[0,\sigma_1]$ belongs to the
$U_{\delta}(id)$-neighborhood of the identity in $\mathrm{Diff}^1_+([0,\sigma_1])$ (the latter group
is being identified with $\mathrm{Diff}^1_+([0,1])$ just by rescaling the interval). Similarly, there is 
$\sigma_2 > 0$ such that $1 - \sigma_2$ is fixed by $H$ and $b$ restricted to $[1-\sigma_2,1]$ 
belongs to the $U_{\delta}(id)$-neighborhood of the identity in $\mathrm{Diff}^1_+([1-\sigma_2,1])$.

For $x \in [0,1]$, let us consider the displacement vector $\triangle(x)$ defined by
\begin{equation}
\triangle (x) := \big( b_1(x)-x, \ldots, b_d(x)-x \big) \in \R^d,
\label{displacement-vector}
\end{equation}
and let us denote by $\| \triangle(x) \|$ its max norm. Notice that
$\| \triangle(x) \| \leq 1$ for all $x\in [0,1]$.

\vspace{0.14cm}

\begin{lem}\label{triconA} {\em For all $r > 0$, there exists $\sigma>0$ such that
$$\triangle(a^{-1}(x))=Da^{-1}(0) \,
A^T \triangle(x)+\epsilon(x) \quad  \text{for all} \quad x\in (0,\sigma)$$
and
$$\triangle(a(x))=Da(1) \,(A^T)^{-1}\triangle(x)+\hat{\epsilon}(x)
\quad \text{for all} \quad  x\in (1-\sigma,1),$$
where \hspace{0.05cm} $\| \epsilon(x) \| \leq r \big( \| \triangle(x) \| + \| \triangle (a^{-1}(x)) \| \big)$ 
\hspace{0.05cm} and \hspace{0.05cm} 
$\|\hat{\epsilon}(x)\| \leq r \big( \| \triangle(x) \| + \| \triangle (a^{-1}(x)) \| \big)$.}
\label{lem:not-so-easy}
\end{lem}

\noindent{\em Proof:} Both assertions being analogous, we will prove only the first one. 
Let $q \in \mathbb{N}$ be such that $\beta_{i,j} := q\alpha_{i,j}$ is an integer for each 
$i,j$. Let $U$ be a neighborhood of the identity in $\mathrm{Diff}^1_+([0,1])$ for which 
Proposition \ref{prop bonatti} holds for $\eta > 0$ and $k \in \mathbb{N}$ defined as 
$$\eta := \frac{r}{2D^2 + \max_i \sum_j | \alpha_{j,i} |}
\quad \mbox{ and } \quad 
k := \max \Big\{ \max_i \Big\{ \sum_j |\beta_{j,i}| \Big\}, q \Big\},$$
where $D := \sup_z \max \{ Da (z), Da^{-1}(z) \}$.
Let $\sigma > 0$ be fixed by $H$ such that the (renormalized) restrictions to $[0,\sigma]$ of 
all the maps $b_j$ and $ab_ja^{-1}$, as well as their inverses, belong to $U$, and such that 
\hspace{0.05cm} $|Da^{-1}(z) - Da^{-1}(0)| \leq \eta$ \hspace{0.05cm} holds for all $z \in [0,\sigma]$. 
Then, by Proposition \ref{prop bonatti},
$$a b_i^q a^{-1} (x) - x = (a b_i a^{-1})^q (x) - x = q ( ab_i a^{-1}(x) - x) + r_{i,1} (x)$$
and
$$b_1^{\beta_{1,i}} \cdots b_d^{\beta_{d,i}} (x) - x = \sum_j \beta_{j,i} (b_j(x) - x) + r_{i,2}(x),$$
where 
%\begin{equation}\label{estimates-r}
$$| r_{i,1} (x) | \leq \eta \big| ab_ia^{-1}(x) - x \big| 
\quad \mbox{ and } \quad 
| r_{i,2}(x) | \leq \eta \max_j \big\{ | b_j(x) - x | \big\} = \eta \| \triangle (x) \|.$$
%\end{equation}
Since
%\begin{equation}\label{conju}
$$a b_i^q a^{-1} = b_1^{\beta_{1,i}} \cdots b_d^{\beta_{d,i}},$$
%\end{equation}
we conclude that
$$q ( ab_i a^{-1}(x) - x) + r_{i,1} (x) = \sum_j \beta_{j,i} (b_j(x) - x) + r_{i,2}(x),$$
hence
\begin{equation}\label{eq linealized} ab_i a^{-1}(x) - x = \sum_j \alpha_{j,i} (b_j(x) - x) + \frac{r_{i,2}(x) - r_{i,1}(x)}{q}.\end{equation}

The $i$-th entry of the vector $\triangle(a^{-1}(x))$ is
$$b_ia^{-1}(x)-a^{-1}(x) = a^{-1}ab_ia^{-1}(x)-a^{-1}(x) = Da^{-1}(z_i) \big( ab_ia^{-1}(x) - x \big),$$
where the last equality holds for a certain point $z_i\in I_x$. By (\ref{eq linealized}) above, for
$x \in (0,\sigma)$, this expression equals
$$Da^{-1}(z_i) \sum_j \alpha_{j,i} \big( b_j(x)-x \big)$$
up to an error $\tilde{\varepsilon}_i(x)$ satisfying
$$|\tilde{\epsilon}_i(x)| \leq Da^{-1} (z_i) \cdot \frac{| r_{i,1} (x) | + | r_{i,2} (x) |}{q} 
\leq 2D \eta \max \big\{ \| \triangle (x) \|, |ab_ia^{-1}(x)-x| \big\}.$$
Since 
$$ab_ia^{-1}(x) - x = a(b_ia^{-1}(x)) - a(a^{-1}(x)) = Da (z_i') \big( b_i a^{-1}(x)) - a^{-1}(x) \big)$$
for a certain $z_i' \in I_{a^{-1}(x)}$, we have
$$|\tilde{\epsilon}_i(x)| 
\leq 
2D \eta \max \big\{ \| \triangle (x) \|, D \| \triangle (a^{-1}(x)) \| \big\} 
\leq 
2D^2 \eta \max \big\{ \| \triangle(x) \|, \| \triangle (a^{-1}(x)) \| \big\}.$$
Moreover, by the choice of $\sigma$, the value of
$Da^{-1}(z_i) \sum_j \alpha_{j,i} (b_j(x)-x)$ equals
$$Da^{-1}(0) \sum_j \alpha_{j,i} (b_j(x)-x)$$
up to an error bounded from above by 
$$\eta \Big| \sum_j \alpha_{j,i} (b_j(x)-x) \Big| 
\leq 
\eta \| \triangle(x) \| \sum_j |\alpha_{j,i}|.$$
Summarizing, $b_ia^{-1}(x)-a^{-1}(x)$ coincides with 
$Da^{-1}(0) \sum_j \alpha_{j,i} \big( b_j(x) - x \big)$
up to an error $\varepsilon_i(x)$ satisfying
$$|\varepsilon_i(x)| 
\leq 2 D^2 \eta \max \big\{ \| \triangle(x) \|, \| \triangle (a^{-1}(x)) \| \big\} 
+ 
\eta \| \triangle(x) \| \sum_j |\alpha_{j,i}|.$$
By the choice of $\eta$, 
%for every $i \in \{1,\ldots,d\}$, 
the last expression is smaller than or equal to $r \big( \| \triangle(x) \| + \| \triangle (a^{-1}(x)) \| \big)$, 
which finishes the proof.
$\hfill\square$

\vs
\vs

Before stating our next lemma, we observe that Lemma \ref{lem tauxconA}
and the compactness of the unit sphere $S^{d-1} \subset \R^d$ imply that for each point
$x_0$ not fixed by $H$, the vectors $\vec{\tau}_{a^{-n}(x_0)}$  (resp., $\vec{\tau}_{a^{n}(x_0)}$)
accumulate at some $\vec{\tau} \in S^{d-1}$ (resp., $\vec{\tau}_*$) as $n \to \infty$. For
each $n \in \mathbb{Z}$, we let $x_n := a^{-n} (x_0)$, and we choose a sequence of
positive integers $n_k$ such that $\vec{\tau}_{x_{n_k}}\to \vec{\tau}$ and
$\vec{\tau}_{x_{-n_k}}\to \vec{\tau}_*$ as $k\to \infty$.

\begin{lem}\label{tauestri} {\em For every $\eta>0$, there exists $K$ such that $k \geq K$ implies
$$\frac{\triangle(x_{n_k})}{\|\triangle(x_{n_k})\|}=\vec{\tau}+\epsilon(k) \quad \mbox{ and }
\quad \frac{\triangle(x_{-n_k})}{\|\triangle(x_{-n_k})\|}=\vec{\tau}_*+\epsilon_* (k),$$
where $\| \epsilon(k) \| \leq \eta $ and $\| \epsilon_*(k) \| \leq \eta$. }
\end{lem}

\noindent{\em Proof:} Let $H_1$ the subgroup of $H$ generated by $b_1,\ldots,b_d$. 
Up to passing to a subsequence if necessary, there is
$b_* \in \{b_1,\ldots, b_d\}$ such that for all $k$,
$$| b_*(x_{n_k})-x_{n_k} | = \max_i \big\{ | b_i(x_{n_k})-x_{n_k} | \big\}.$$
Then the functions $\psi_k: H_1 \to \R$ defined by
$$\psi_k (b) = \frac{b(x_{n_k}) - x_{n_k}}{b_*(x_{n_k})-x_{n_k}}$$
converge as $k \to \infty$ to a homomorphism $\psi \!: H_1 \to \R$ which is normalized, in the
sense that $\max_i |\psi(b_i)| = 1$. Indeed, this is the content of the Thurston's stability theorem
\cite{thurston} (which in its turn can be easily deduced from Proposition \ref{prop bonatti}).

The vectors $\vec{\tau}_k$ and $\vec{\tau}$ naturally induce normalized homomorphisms
from $H$ into $\mathbb{R}$, namely the {\em normalized} translation number homomorphisms, 
and their {\em limit} homomorphism. We denote them by $\vec{\tau}_k$ and $\vec{\tau}$, respectively. 
For these homomorphisms and any $b,c$ in $H$, the inequality $\vec{\tau} (b) < \vec{\tau} (c)$
implies $\vec{\tau}_k (b) < \vec{\tau}_k (c)$ for $k$ larger than a certain $K_0$,
which implies $b (z) < c (z)$ for all
$z \in I_{x_{n_{k}}}$ and all $k > K_0$. By evaluating at $z = x_{n_k}$,
this yields $\psi_k (b) < \psi_k (c)$ for $k > K_0$.
Passing to the limit, we finally obtain $\psi(b) \leq \psi(c)$.
As a consequence, there must exist a constant $\lambda$ for
which $\vec{\tau} = \lambda \psi$. Nevertheless, since both
homomorphisms are normalized (and point in the same
direction), we must have $\lambda = 1$, which yields
the convergence of $\triangle (x_{n_k}) / \| \triangle (x_{n_k})\|$
towards $\vec{\tau}$. The second convergence is proved in an analogous way. $\hfill\square$

\vs
\vs

Henceforth, and in many other parts of this work, we will use a trick due to
Muller and Tsuboi that allows reducing to the
case where all group elements are tangent to the identity at the endpoints. 
This is achieved after conjugacy by an appropriate homeomorphism that is smooth 
at the interior and makes  flat the germs at the enpoints. In concrete terms, we have:

\begin{lem} [\cite{muller,tsuboi}]
{\em Let us consider the germ at the origin of the local (non-differentiable) homeomorphism 
$\varphi (x) := sgn(x) \exp(-1/|x|)$. Then for every germ of $C^1$ diffeomorphism $f$ (resp. 
vector field $\mathcal{X}$) at the origin, the germ of the conjugate $\varphi^{-1} \circ f \circ \varphi$ 
(resp., push-forward $\varphi_* (\mathcal{X})$) is differentiable and flat in a neighborhood 
of the origin.} \label{MT-trick}
\end{lem}

We should stress, however, that although this lemma simplifies many computations, 
in what follows it may avoided just noticing that, as $Da$ is continuous,  the element 
$a$ behaves like an affine map close to each endpoint.

\vsp

Recall that $\R^d$ decomposes as $E^s\oplus E^u$, where $E^s$ (resp.
$E^u$) stands for the stable (unstable) subspace of $A^T$. We denote by
$\pi_s$ and $\pi_u$ the projections onto $E^s$ and $E^u$, respectively.
We let $\| \cdot \|_*$ be the natural norm on $\R^d$ associated to this
direct-sum structure, namely,
$$\| v \|_* := \max \big\{\| \pi_s (v) \|, \|\pi_u (v) \| \big\}.$$

\begin{lem} \label{lem proyectivo}{\em
For any neighborhood $V \subset S^{d-1}$ of $E^u \cap S^{d-1}_*$ in the unit sphere
$S^{d-1}_* \subset \R^d$ (with the norm $\|\cdot\|_*$), there is $\sigma>0$
such that for all $x \in (0,\sigma)$ not fixed by $H$,
$$\frac{\triangle(x)}{\|\triangle(x)\|_*}\in V \large\implies
\frac{\triangle(a^{-1}(x))}{\| \triangle(a^{-1}(x)) \|_*} \in V.$$
Moreover, if $V$ is small enough, then there exists $\kappa > 1$ such that
$$\frac{\triangle(x)}{\| \triangle(x) \|_*} \in V
\large\implies
\| \triangle(a^{-1}x) \|_*\geq
 \kappa \|\triangle(x)\|_* .$$}
\end{lem}

\noindent{\em Proof:} For the first statement, we need to show that for every prescribed
positive $\varepsilon < 1$, for points $x$ close to the origin and not fixed by $H$, we have
$$\frac{\|\pi_s \triangle (a^{-1}(x))\|}{\|\pi_u \triangle (a^{-1}(x)) \|} < \varepsilon
\quad \mbox{ provided } \quad
\frac{\| \pi_s \triangle (x) \|}{\| \pi_u \triangle (x) \|} < \varepsilon.$$

Let $\lambda > 1$ (resp., $\lambda' < 1$) be such that the norm of nonzero vectors
in $E^u$ (resp., $E^s$) are expanded by at least $\lambda$ (resp., contracted by at least
$\lambda'$) under the action of $A^T$. Choose a positive $r < \varepsilon$ small enough so that
\begin{equation}\label{cond-on-r}
\frac{1 + r}{1 - r} 
\left[ \frac{\varepsilon \lambda'}{\lambda - 2r} + \frac{2r}{\lambda - 2r} \right] 
\leq 
\frac{\varepsilon - \frac{r}{1 - r}}{1 + \frac{\varepsilon r}{1 + r}}.
\end{equation}
Consider a point $x$ not fixed by $H$ lying
in the interval $(0,\sigma)$ given by Lemma \ref{lem:not-so-easy}. Then from
\begin{eqnarray*}
\| \pi_s \triangle (a^{-1}(x)) \| 
&\leq& 
\| \pi_s A^T \triangle (x)\| + r \big( \| \triangle(x) \| + \| \triangle (a^{-1}(x)) \| \big) \\
&\leq& 
\lambda' \| \pi_s \triangle(x) \| + r \| \pi_s \triangle(x)\| + r \| \pi_u \triangle(x)\| + r \| \pi_s \triangle(a^{-1}(x))\| + r \| \pi_u \triangle(a^{-1}(x))\| 
\end{eqnarray*}
we obtain
\begin{equation}\label{first-comp}
\| \pi_s \triangle (a^{-1}(x)) \| 
\leq 
\frac{1}{1-r} \left[ \lambda' \| \pi_s \triangle(x) \| + 2r \| \pi_u \triangle(x)\| + r \| \pi_u \triangle(a^{-1}(x))\| \right].
\end{equation}
Similarly, from 
\begin{eqnarray*}
\| \pi_u \triangle (a^{-1}(x)) \| 
&\geq& \| \pi_u A^T \triangle (x)\| - r \big( \| \triangle(x) \| + \| \triangle (a^{-1}(x)) \| \big) \\
&\geq& \lambda \| \pi_u \triangle(x) \| -  2r \| \pi_u \triangle(x) \| - r \| \pi_u \triangle (a^{-1}(x)) \| - r \| \pi_s \triangle (a^{-1}(x)) \|, \\
\end{eqnarray*}
we obtain
\begin{equation}\label{second-comp}
\| \pi_u \triangle (a^{-1}(x)) \| 
\geq
\frac{1}{1+r} \left[ (\lambda - 2r) \| \pi_u \triangle (x) \| - r \| \pi_s \triangle (a^{-1}(x)) \|  \right].
\end{equation}
Thus, letting $\alpha := \| \pi_s \triangle (a^{-1}(x)) \| $, $\beta := \| \pi_u \triangle (a^{-1}(x)) \| $, 
$$A := \frac{1}{1-r} \left[ \lambda' \| \pi_s \triangle(x) \| + 2r \| \pi_u \triangle(x)\| \right] 
\quad \mbox{ and } \quad 
B := \frac{(\lambda - 2r) \| \pi_u \triangle (x) \|}{1+r},$$
we have that (\ref{first-comp}) and (\ref{second-comp}) translate into 
$$\alpha -  \frac{\beta r}{1-r} \leq A \qquad \mbox{ and } \qquad \beta + \frac{\alpha r}{1 + r} \geq B,$$
and hence,
$$\frac{\frac{\alpha}{\beta} - \frac{r}{1-r}}{1 + \frac{\alpha}{\beta} \cdot \frac{r}{1+r}} \leq \frac{A}{B}.$$
%and hence $\frac{\alpha}{\beta}-\frac{r}{1-r}\leq \frac{A}{B}+\frac{\alpha}{\beta}\frac{A}{B}\frac{r}{1+r}$. From where we obtain that
%which easily implies that 
%$$\frac{\alpha}{\beta} \leq \frac{\frac{A}{B} + \frac{r}{1-r}}{1 - \frac{A}{B} \cdot \frac{r}{1+r}}.$$
But 
$$\frac{A}{B} 
= 
\frac{1+r}{1-r} \left[ 
\frac{\lambda' }{\lambda - 2r} \frac{\| \pi_s \triangle (x) \|}{ \| \pi_u \triangle (x) \| } + \frac{2r}{\lambda - 2r} \right] 
<
\frac{1+r}{1-r} \left[ 
\frac{\varepsilon \lambda' }{\lambda - 2r} + \frac{2r}{\lambda - 2r} \right].$$
Therefore, by the  choice of $r$ (see (\ref{cond-on-r})), and the fact that $ x\mapsto \frac{x-c}{1+c d}$ is an increasing on $x$ 
for positive $c,d$, we obtain that $\alpha / \beta < \varepsilon$, 
which shows the first assertion of the lemma.

To conclude the proof, notice that by the estimate (\ref{second-comp}) above,
$$ \Big( 1 + \frac{r}{1+r} \Big) \| \triangle (a^{-1}(x)) \|_* 
=
\Big( 1 + \frac{r}{1+r} \Big) \| \pi_u \triangle (a^{-1}(x)) \| 
\geq 
\frac{\lambda - 2r}{1 + r} \| \pi_u \triangle(x) \| 
= 
\frac{\lambda - 2 r}{1 + r} \| \triangle (x) \|_*,$$
which shows the second assertion of the lemma for $\kappa := (\lambda - 2r) / (1 + 2 r)$ and $r$ small enough. 
$\hfill\square$

\vs
\vs

Now we can easily finish the proof of Proposition \ref{prop no hay
accione por nivels}. To do this, choose a point $x_0\in (0,1)$ that
is not fixed by $H$. We need to consider two cases:

\vs

\noindent {\bf Case 1:} $\vec{\tau}_{x_0} \notin E^s$

\vsp

In this case, we first observe that Lemma \ref{lem:vary} implies that any accumulation
point of $\vec{\tau}_{a^n(x_0)}$ (in particular, $\vec{\tau}$) must belong to $E^u$. Let 
$V$ be a small neighborhood around $E^u\cap S^{d-1}_*$ in $S^{d-1}_*$ so that 
both statements of Lemma \ref{lem proyectivo} hold. Then, by Lemma \ref{tauestri},
the vector $\triangle (x_k) / \| \triangle (x_k) \|_*$ belongs to $V$ starting
from a certain $k = K$. This allows applying Lemma \ref{lem proyectivo}
inductively, thus showing that for all $n \geq 0$,
$$1\geq \| \triangle(x_{n+k}) \|_* \geq \kappa^n \| \triangle (x_k) \|_*.$$
Letting $n$ go to infinity, this yields a contradiction.

\vs

\noindent {\bf Case 2:} $\vec{\tau}_{x_0}\in E^s$.
\vs

In this case, Lemma \ref{lem:vary} yields $\vec{\tau}_* \in E^s$. We then
proceed as above but on a neighborhood of 1 working with $a$
instead of $a^{-1}$ and with $(A^T)^{-1}$ instead of $A^T$. Details
are left to the reader. (This requires for instance an analog of Lemma
\ref{lem proyectivo} for the dynamics close to $1$.)

%%%%%%%%%%%%%%%%%%%%%%%%%%%%%%%%%%%%%%%%%%%%%%%%%

\subsection{Minimality of affine-like actions}

%%%%%%%%%%%%%%%%%%%%%%%%%%%%%%%%%%%%%%%%%%%%%%%%%

In this section, we begin by showing Proposition \ref{prop conjugacion}.
Let $\phi: G \to \mathrm{Diff}^1_+([0,1])$ be a representation
with non-Abelian image. We know from Proposition
\ref{prop no hay accione por nivels} that $\phi$ is semiconjugate 
to a representation $\psi \!: G \to \mathrm{Aff}_+([0,1])$ in
the affine group. The elements in the commutator subgroup $[\psi(G),\psi(G)]$
are translations.
In what follows, we will assume that the right endpoint is topologically attracting
for $\psi(a)$, hence $\psi(a)$ is conjugate to an homothety $x\to \lambda x$
with $\lambda>1$ (the other case is analogous).
Up to changing $a$ by a positive power, we may assume
that $\lambda \geq 2$. We fix $b \!\in\! H$ such that $\psi(b)$ is a non-trivial
translation. Up to changing $b$ by its inverse and conjugating $\psi$ by
an appropriate homothety, we may assume that $\psi(b) = T_1$.
We consider a finite system of generators of $G$
that contains both $a$ and $b$.

Suppose for a contradiction that $\phi(G)$ does not act minimally. Then there
is an interval $I$ that is wandering for the action of $[\phi(G),\phi(G)]$.
As before, we may assume that $D\phi(c)(1)=1$ for all $c\in G$.
Fix $\varepsilon > 0$ such that $ (1-\varepsilon)^3 > 1/2$, 
and let $\delta>0$ be such that
\begin{equation}
\label{eq derivada acotada}1-\epsilon \leq D \phi(c) (x) \leq 1+\epsilon
\quad \mbox{for each } c \in \{ a^{\pm 1}, b \} \mbox{ and all } x\in [1-\delta, 1].
\end{equation}
Clearly, we may assume that $I\subset [1-\delta,1]$.

Notice that $\psi(a^{-k} b a^k) = T_{\lambda^{-k}}$ for all $k \in \mathbb{Z}$.
We consider the following family of translations
$$h_{(\varepsilon_i)} := (a^{-n}b^{\varepsilon_n} a^{n}) \cdots
(a^{-2} b^{\varepsilon_2} a^{2}) (a^{-1} b^{\varepsilon_1} a),$$
%= a^{-n} (h^{\varepsilon_n}a) \cdots (h^{\varepsilon_1} a),$$
where $(\varepsilon_i) = (\varepsilon_1, \ldots,\varepsilon_n) \in \{0,1\}^n$.
These satisfy the following properties:

\vspace{0.25cm}

\noindent (i) We have that
$(\varepsilon_i) \not= (\tilde{\varepsilon}_i)$ implies
$h_{(\varepsilon_i)} \not = h_{(\tilde{\varepsilon}_i)}$:
this easily follows from that $\lambda \geq 2$.

\vspace{0.25cm}

\noindent (ii) We have $\phi(h_{(\varepsilon_i)})(1-\delta) \geq 1-\delta$: this follows
from that $\phi(b)$ attracts towards $1$ and that $\varepsilon_i \geq 0$ for all $i$.

\vspace{0.25cm}

\noindent (iii) The element \hspace{0.01cm} 
$h_{(\varepsilon_i)} = a^{-n} (b^{\varepsilon_n}a)
\cdots (b^{\varepsilon_2} a) (b^{\varepsilon_1} a)$ 
\hspace{0.01cm} belongs to the ball of radius
$3n$ in $G$. In particular, due to (\ref{eq derivada acotada}) and the
preceding claim, we have \hspace{0.01cm}
$D \phi(b_{(\varepsilon_i)})(x) \geq (1-\epsilon)^{3n}$ 
\hspace{0.01cm} for all $x\in [1-\delta,1]$.

\vspace{0.25cm}

Since for each $c \in G$ there exists $x_I \in I$ for which
$|c(I)| = Dc(x_I)|I|$ (where $| \cdot  |$ stands for the length of the
corresponding interval), putting together the assertions above we conclude
$$1 \geq \sum_{(\varepsilon_i)} \big| h_{(\varepsilon_i)} (I) \big|
\geq 2^n (1-\epsilon)^{3n} |I| > 1,$$
where the last inequality holds for $n$ large enough. This
contradiction finishes the proof of Proposition \ref{prop conjugacion}.

\vspace{0.3cm}

It should be emphasized that Proposition \ref{prop conjugacion} is no longer true 
for $C^1$ (even real-analytic) actions on the real line (equivalently, on the open 
interval). Indeed, this issue was indirectly adressed by Ghys and Sergiescu in 
\cite[section III]{ghys sergiescu}, as we next state and explain.

\vspace{0.1cm}

\begin{thm}\label{teo semiconj} (\cite{ghys sergiescu}).
{\em The Baumslag-Solitar group $BS(1,2):= \big\langle a,b\mid aba^{-1}=b^2 \big\rangle$ 
embeds into $\mathrm{Diff}_+^1(\R)$ via an action that is semiconjugate, but not conjugate, 
to the canonical affine action and such that the element $a \!\in\! B(1,2)$ acts with two fixed points.}
\end{thm}

\vspace{0.1cm}

Recall that $BS(1,2)$ is isomorphic to the group of order-preserving
affine bijections of $\Q_2$, where $\Q_2$ is the group of diadic rationals.
Hence, every element in $BS(1,2)$ may be though as a pair
$\big( 2^n,\frac{p}{2^q} \big)$, which identifies to the affine map
$$\left( 2^n,\frac{p}{2^q} \right) \!: \hspace{0.1cm} x\to 2^nx+\frac{p}{2^q}.$$
Notice that $\Q_2$ corresponds to the subgroup of translations inside $BS(1,2)$.

Next, following \cite{ghys sergiescu}, we consider a homeomorphism
$f:\R\to \R$ satisfying the following properties:

\vspace{0.1cm}

\noindent (I) For every $x\in \R$, we have $f(x+1)=f(x)+2$.

\vspace{0.1cm}

\noindent (II) $f(0)=0$.

\vspace{0.1cm}

\begin{lem}[\cite{ghys sergiescu}] {\em The map $\theta_f:\frac{p}{2^q}
\in \Q_2\to f^{-q}T_pf^q \in \mathrm{Homeo}_+(\R)$ is a
well-defined homomorphism.}

\end{lem}

\begin{lem}[\cite{ghys sergiescu}] The map
$\big( 2^n,\frac{p}{2^q} \big) \in BS(1,2) \to \theta_f(\frac{p}{2^q})
\circ f^n \in \mathrm{Homeo}_+(\R)$ is a group homomorphism.
%This homomorphism will be still denoted by $\theta_f$.
\end{lem}

The homomorphism provided by the last lemma above will still be denoted by $\theta_f$.
Notice that $\theta_f(a)=f$.

Next, for $1\leq r\leq \infty, \omega$, we impose a third condition on $f$:

\vspace{0.1cm}

\noindent (III$_r$) The map $f$  is of class $C^r$.

\vspace{0.1cm}

We have

\begin{lem}
[\cite{ghys sergiescu}] The image $\theta_f(BS(1,2))$ is a subgroup of
$\mathrm{Diff}^r_+(\R)$.
\end{lem}

We end with

\begin{lem}[\cite{ghys sergiescu}] Suppose that the function $f$ has at least two fixed
points. Then $\theta_f(BS(1,2))$ has an exceptional minimal set ({\em i.e.} a minimal
invariant closed set locally homeomorphic to the Cantor set).
\end{lem}

To close this section, we point out that a similar construction
can be carried out for all Baumslag-Solitar's groups
$BS(1,n) := \big\langle a,b \mid bab^{-1} = a^n \big\rangle$.
Roughly, we just need to replace condition (I) by:

\vsp

\noindent (I)$_n$ For every $x \in \R$, we have $f(x+1) = f(x+n)$.

%%%%%%%%%%%%%%%%%%%%%%%%%%%%%%%%%%%%%%%%%%%%%%%%%%%%%%%%%%%%%%%%

\subsection{Rigidity of multipliers}

We start by dealing with the Baumslag-Solitar group $BS(1,2)$. Let us consider a faithful
action of this group by $C^1$ diffeomorphisms of the closed interval. We known
that such an action must be topologically conjugate to an affine action, hence to the
standard affine action given by $a \!: x \mapsto 2x$ and $b \! : x \mapsto x+1$. (It
is not hard to check that all faithful affine actions of $B(1,2)$ are conjugate inside  
$\mathrm{Aff}(\R)$.) Let $\varphi: (0,1) \rightarrow \R$ denote the topological 
conjugacy. Our goal is to show

\vsp

\begin{prop} \label{derivada=2}
The derivative of $a$ at the interior fixed point equals 2.
\end{prop}

\noindent{\em Proof:} For the proof, we let $I := \varphi^{-1} ([0,1])$.
Notice that for all positive integers $n,N$, the intervals
$$(a^{-n} b^{\varepsilon_n} a^n) \cdots (a^{-2} b^{\varepsilon_2} a^2)
(a^{-1}b^{\varepsilon_1}a) b^N a^{-n}(I), \quad \varepsilon_i \in \{0,1\},$$
have pairwise disjoint interiors. Indeed, these intervals are mapped by
$\varphi$ into the dyadic intervals of length $1/2^n$ contained in $[N,N+1]$.
For simplicity, we assume below that both $a$ and $b$ have derivative 1 at
the endpoints. (As before, this may be performed via the Muller-Tsuboi trick; 
{\em c.f.} Lemma \ref{MT-trick}). 

Assume first that $Da (x_0) < 2$, where $x_0$ is the interior fixed point of $a$.
Then there are $C > 0$ and $\varepsilon > 0$ such that for all $n \geq 1$,
$$\big| a^{-n}(I) \big| \geq C \Big( \frac12 + \varepsilon \Big)^n |I|.$$
Fix $\delta > 0$ such that
\begin{equation}\label{unila}
\big( 1-\delta \big)^3 \Big( \frac12 + \varepsilon \Big) > 1/2.
\end{equation}
Let $\sigma > 0$ be small enough so that
$$Da(x) \geq 1-\delta, \quad Da^{-1}(x) \geq 1-\delta  \quad \mbox{and} \quad Db(x) \geq 1-\delta
\quad \mbox{for all} \quad x \in [1-\sigma,1].$$
Finally, let $N \geq 1$ be such that $b^N(x_0) \geq 1 - \sigma$. Similarly
to the proof of Proposition \ref{prop conjugacion}, for such $N$ and all
$n \geq 1$, we have for all choices $\varepsilon_i \in \{0,1\}$,
$$\big| (a^{-n} b^{\varepsilon_n} a^{n}) \cdots (a^{-2} b^{\varepsilon_2} a^{2})
(a^{-1}b^{\varepsilon_1}a) b^N a^{-n}(I) \big|
\geq (1-\delta)^{3n} DC\Big(\frac12 +\varepsilon \Big)^n \big| I \big|,$$
where $D := \min_x Db^N (x)$. As there are $2^n$ of these intervals, we have
$$1 \geq 2^n (1-\delta)^{3n} DC \Big( \frac12+\varepsilon \Big)^n |I|,$$
which is impossible for a large-enough $n$ due to (\ref{unila}).

Assume next that $Da (x_0) > 2$. Then there are $C' > 0$ and $\varepsilon' > 0$
such that for all $n \geq 1$,
$$\big| a^{-n}(I) \big| \leq C' \Big( \frac12 - \varepsilon' \Big)^n.$$
Fix $\delta' > 0$ such that
\begin{equation}\label{dorila}
\big( 1+\delta' \big)^3 \Big( \frac12 - \varepsilon' \Big) < 1/2.
\end{equation}
Let $\sigma' > 0$ be small enough so that
$$Da(x) \leq 1+\delta,   \quad Da^{-1}(x) \leq 1+\delta  \quad\mbox{and} \quad Db(x) \leq 1+\delta
\quad \mbox{for all} \quad x \in [1-\sigma',1].$$
Finally, let $N' \geq 1$ be such that $b^{N'}(x_0) \geq 1 - \sigma'$.
Proceeding as before, we see that for such $N'$ and all $n \geq 1$,
we have for all choices $\varepsilon_i \in \{0,1\}$,
$$\big| (a^{-n} b^{\varepsilon_n} a^{n}) \cdots (a^{-2} b^{\varepsilon_2} a^{2})
(a^{-1}b^{\varepsilon_1}a) b^N a^{-n}(I) \big|
\leq (1+\delta')^{3n} D'C'\Big(\frac12 -\varepsilon' \Big)^n \big| I \big|,$$
where $D' := \max_x Db^N (x)$. However, the involved intervals cover
$I_{N'} := b^{N'} (I) = \varphi^{-1} \big( [N',N'+1] \big)$. Thus,
$$ |I_{N'}|
\leq 2^n (1+\delta')^{3n} D'C'\Big(\frac12 -\varepsilon' \Big)^n \big| I \big|,$$
which is again impossible for a large-enough $n$ due to $(\ref{dorila})$.
$\hfill\square$

\vsp
\vsp

\begin{rem} The action of the Baumslag-Solitar group by $C^1$ diffeomorphisms 
of the real line constructed in the preceding section can be easily modified into a 
minimal one for which the derivative of $a$ at the fixed point equals 1. Roughly, 
we just need to ask for the map $f$ along the construction to have a single fixed 
point, with derivative 1 at this point. This shows that Theorem~\ref{thm multiplier} 
is no longer true for actions by $C^1$ diffeomorphisms of the open interval.
\end{rem}

\vsp

The preceding proposition corresponds to a particular case of Theorem
\ref{thm multiplier} but illustrates the technique pretty well. Below we give
the proof of the general case along the same ideas. First, as $A$ is supposed to be
hyperbolic, we know that the action of $G$ is topologically conjugate to an affine
one. Moreover, Proposition \ref{propafin} completely describes such an action:
up to a topological conjugacy $\varphi$, it is given by
correspondences $a \mapsto M_\lambda$ and $h_i \mapsto T_{t_i}$,
where $(t_1,\ldots,t_d)$ is an eigenvector of $A$ with eigenvalue $\lambda$.
Up to conjugacy in $\mathrm{Aff} (\R)$, we may assume that one of the
$t_i's$ equals 1, hence $b:= b_i$ is sent into $T_t := T_1$.

Next, we proceed as above, but with a little care. Notice that changing $a$ by
an integer power if necessary, we may assume that $\lambda \geq 2$.

Assume first that
$Da (x_0) < \lambda$, where $x_0$ is the interior fixed point of $a$.
Then there are $C > 0$ and $\varepsilon > 0$ such that for all $n \geq 1$,
$$\big| a^{-n}(I) \big| \geq C \Big( \frac{1}{\lambda} + \varepsilon \Big)^n.$$
Fix $\delta > 0$ such that \hspace{0.01cm} $(1-\delta)^3 (\frac{1}{\lambda} + \varepsilon)
> \frac{1}{\lambda}$. \hspace{0.01cm} Let $\sigma > 0$ be small so that 
$Da(x) \geq 1-\delta$, $Da^{-1}(x) \geq 1-\delta \mbox{ and } Db(x) \geq 1-\delta
\mbox{ hold for all }  x \in [1-\sigma,1].$
Finally, let $N \geq 1$ be such that $b^N(x_0) \geq 1 - \sigma$. Given
$n \geq 1$, we consider for all choices $\varepsilon_i \in \{0,1,\ldots,[\lambda]\}$,
the intervals $(a^{-n} b^{\varepsilon_n} a^{n}) \cdots (a^{-2} b^{\varepsilon_2} a^{2})
(a^{-1}b^{\varepsilon_1}a) b^N a^{-n}(I)$, where $I$ is the preimage of $[0,1]$ under
the topological conjugacy into the affine action.  As before, we have for each such choice
$$\big| (a^{-n} b^{\varepsilon_n} a^{n}) \cdots (a^{-2} b^{\varepsilon_2} a^{2})
(a^{-1}b^{\varepsilon_1}a) b^N a^{-n}(I) \big|
\geq (1-\delta)^{3n} DC\Big(\frac{1}{\lambda} + \varepsilon \Big)^n \big| I \big|,$$
where $D := \min_x Db^N (x)$. These intervals do not necessarily have pairwise 
disjoint interiors, but their union covers $I$ with multiplicity at most 2. As 
there are $\big( [\lambda]+1 \big)^n$ of these intervals, we have
$$2 \geq \big( [\lambda]+1 \big)^n (1-\delta)^{3n}
DC \Big( \frac{1}{\lambda}+\varepsilon \Big)^n |I|,$$
which is impossible for large-enough $n$.

Assume next that $Da (x_0) > \lambda$.
Then there are $C '> 0$ and $\varepsilon' > 0$ such that for all $n \geq 1$,
$$\big| a^{-n}(I) \big| \leq C' \Big( \frac{1}{\lambda} - \varepsilon' \Big)^n.$$
Fix $\delta' > 0$ such that $(1+ \delta')^3 (\frac{1}{\lambda} - \varepsilon')
\!<\! \frac{1}{\lambda}$. Let $\sigma' > 0$ be small enough so that
$Da(x) \leq 1+\delta'$, $Da^{-1}(x) \leq 1+\delta'$ and 
$Db(x) \!\leq \!1+\delta' \mbox{ hold for all } x \in [1-\sigma',1].$
Finally, let $N' \geq 1$ be such that $b^{N'}(x_0) \geq 1 - \sigma'$.
As before, given $n \geq 1$, for all choices
$\varepsilon_i \in \{0,1,\ldots,[\lambda]\}$, we have
$$\big| (a^{-n} b^{\varepsilon_n} a^{n}) \cdots (a^{-2} b^{\varepsilon_2} a^{2})
(a^{-1}b^{\varepsilon_1}a) b^{N'} a^{-n}(I) \big|
\leq (1+\delta')^{3n} D'C'\Big(\frac{1}{\lambda} - \varepsilon' \Big)^n \big| I \big|,$$
where $D' := \max_x Db^N (x)$. These intervals cover $I_{N'} := b^{N'}(I)$ for each
$n \geq 1$. As there are $([\lambda]+1)^n$
of these intervals, we have
$$|I_{N'}|
\leq
\big( [\lambda]+1 \big)^n (1+\delta')^{3n} DC \Big( \frac{1}{\lambda}-\varepsilon' \Big)^n |I|.$$
Although this is not enough to conclude, we notice that we may replace $a$ by
$a^k$ along the preceding computations, now yielding
$$|I_{N'}|
\leq \big( [\lambda^k]+1 \big)^n
(1+\delta')^{3n} DC \Big( \frac{1}{\lambda}-\varepsilon' \Big)^{kn} |I|.$$
Choosing $k$ large enough so that
$$\big( [\lambda^k] + 1 \big) \Big( \frac{1}{\lambda} - \varepsilon' \Big)^k
(1 + \delta')^3 < 1$$
and then letting $n$ go to infinity, this gives the desired contradiction.

We have hence showed that $Da (x_0) = \lambda$. To show that the derivative of
$a^k b$ at the interior fixed point equals $\lambda^k$ for each $k \neq 0$ and
all $b \in H$, just notice that the associated affine action can be conjugate in
$\mathrm{Aff}(\R)$ so that $a^k b$ is mapped into $T_{\lambda^k}$. 
Knowing this, we may proceed in the very same way as above.

%%%%%%%%%%%%%%%%%%%%%%%%%%%%%%%%%%%%%%%%%%%%%%%%%%%%%%%%%%%%%%%%

\subsection{On the smoothness of conjugacies}

As we announced in the Introduction, actions by $C^1$ diffeomorphisms are rarely
rigid in what concerns the regularity of conjugacies. In our context, this is actually
never the case, as it is shown by the next

\vsp

\begin{prop} \label{p:bad-conjugates}
{\em Let $G$ be a finitely-generated group of the form $\Z \ltimes_A H$, where 
$A \in GL_d(\Q)$ and $rank_{\Q}(H) = d$. Then every faithful action of $G$ by $C^1$ diffeomorphisms 
of $[0,1]$ can be approximated in the $C^1$ topology by actions by $C^1$ diffeomorphisms 
that are topologically conjugate to it but for which no Lipschitz conjugacy exists.}
\end{prop}

\vsp

The proof of this proposition follows by an standard application of the Anosov-Katok 
method. The reader is referred to \cite{FH,FK} for a general panorama on this, yet 
the construction we give below is self-contained. It is to be noticed that the only dynamical properties we need are that we 
deal with nontrivial actions for which there are free orbits and no homeomorphism fixing two interior points can conjugate 
the action to itself.
%\marginpar{No entiendo esta frase. Es falsa para la identidad}. 
Thus, the proposition applies in many more situations.

We start with an elementary lemma.

\vsp

\begin{lem} \label{l:familias}
{\em There exists a family of $C^{\infty}$ diffeomorphisms $\varphi_{I}^D \!: I \to I$ between closed intervals $I$, 
where $D \geq 1$, that are infinitely tangent to the identity at the endpoints, satisfy $\varphi_{\ell}^D (m) = m$ and 
$D \varphi_{\ell}^D (m) = D$ for the midpoint $m$ of $I$, and such that given $\varepsilon > 0$ and $\bar{D} > 1$, 
there exists $\delta > 0$ such that for all $y \in I$ and all $D,D'$ in $[1,\bar{D}]$  
%and all $x,y$ in $I$ 
satisfying that $1 - \delta \leq D/D' \leq 1+\delta$, we have }
%and $|x-y| \leq \delta$,} 
$$1 - \varepsilon \leq \frac{D \varphi_{I}^D (y)}{D \varphi_{I}^{D'}(y)} \leq 1 + \varepsilon.$$
\end{lem}

\noindent{\em Proof:} Let $\chi$ be a vector field on $[0,1]$ which is infinitely flat at the endpoints, strictly negative 
(resp. positive) on $[-1,0]$ (resp. $[0,1]$), and satisfies $\chi (x) = x \frac{\partial }{ \partial x}$ in a neighborhood 
of the origin. Then let $\varphi_{[-1,1]}^D$ be the flow associated to $\chi$ up to time $T := \log(D)$, so that 
$D\varphi_{[-1,1]}^D (0) = D$. Finally, for an interval $I$ as in the statement, define $\varphi_{I}^D$ to be the affine conjugate of 
$\varphi_{[-1,1]}^D$. Then the family $\{ \varphi_{I}^D \}$ satisfies all the desired properties, as the reader may easily check. 
$\hfill\square$

\vsp\vsp\vsp\vsp

\noindent{\em Proof of Proposition \ref{p:bad-conjugates}.} Start with $G$ viewed as a group of $C^1$ diffeomorphisms of $[0,1]$, 
and remind that the action of $G$ at the interior is topologically conjugate to a non-Abelian affine group. 
%\marginpar{No decimos que es conjugado a afin? talvez decirlo arriba?}. 
Fix two points $x_0, y_0$ in $(0,1)$ belonging to two different orbits, with $x_0$ having a free orbit by $G$, and fix also a finite generating set of $G$. 
%and denote by $a$ the generator of the $\Z$-factor of $G$ and by $\{b_1,\ldots,b_d\}$ a $\Q$-basis of $\Q^d$, 
%which we assume to be the canonical one and contained in $H$ (see the beginning of \S \ref{section-remark-base-canonica}).
We will inductively construct a sequence of diffeomorphisms $\varphi_k$ of $[0,1]$ in such a way 
%\in \{a^{\pm 1}, b_i^{\pm 1} \}$, where $i \in \{1,\ldots,d\}$,
that we have for $\tilde{\varphi}_k := \varphi_1 \circ \cdots \circ \varphi_k$:

\vsp\vsp

\noindent (i) \hspace{0.35cm}
$\big\| \tilde{\varphi}_{k+1} - \tilde{\varphi}_k \big\|_{C^0} \leq \frac{1}{2^k}$,

\vsp\vsp

\noindent (ii) \hspace{0.22cm}
$\big\| \tilde{\varphi}_{k+1} \circ c \circ \tilde{\varphi}_{k+1}^{-1} -
\tilde{\varphi}_{k} \circ c \circ \tilde{\varphi}_{k}^{-1} \big\|_{C^1} \leq \frac{1}{2^k}$ for each generator $c$,

\vsp\vsp

\noindent (iii) \hspace{0.18cm} $x_0,y_0$ are both fixed by $\varphi_k$,  

\vsp\vsp

\noindent (iv) \hspace{0.12cm}
$D \varphi_{k} (x_0) >  k / \min_y D \tilde{\varphi}_{k-1} (y)$, 
and if we denote by $J_{k}$ the connected component of the set 
$\big\{x \mid D \varphi_{k} (x) > k / \min_y D \tilde{\varphi}_{k-1} (y) \big\}$ 
containing $x_0$, then the support of $\varphi_{k+1}$ has measure $< |J_k|/2$.

\vsp\vsp\vsp

Assume for a while we have performed this construction, and let us complete the proof. 
By (i), we have that the sequence $(\tilde{\varphi}_k)$ converges to a homeomorphism
$\tilde{\varphi}_{\infty}$. By (ii), the sequence of the actions conjugated by $\tilde{\varphi}_k$
converge in the $C^1$ topology to the action conjugated by $\tilde{\varphi}_{\infty}$. Due 
to (iii), each $\tilde{\varphi}_k$ fixes $x_0$ and $y_0$, hence the same holds for
$\tilde{\varphi}_{\infty}$. As conjugacies to affine groups with dense translation
subgroup are unique up to right composition with an affine map, we deduce
that $\tilde{\varphi}_{\infty}$ is the unique conjugacy between $G$ and \hspace{0.04cm}
$\tilde{\varphi}_{\infty} G \tilde{\varphi}_{\infty}^{-1}$ \hspace{0.04cm} fixing these two
points. Finally,  using (iv), it readily follows that the derivative of $\tilde{\varphi}_k$ 
is larger than $k$ 
%\marginpar{Porque desparece el minimo?} 
on certain intervals that remain disjoint from the supports of 
$\varphi_{k+1},\varphi_{k+2}, \ldots$ 
As a consequence, the limit homeomorphism $\tilde{\varphi}_{\infty}$ is not 
Lipschitz. Because of the uniqueness of the conjugacy up to affine transformations
previously discussed, this implies that $G$ and \hspace{0.03cm}
$\tilde{\varphi}_{\infty} G \tilde{\varphi}_{\infty}^{-1}$ \hspace{0.03cm}
cannot be conjugated by any Lipschitz homeomorphism.

To conclude the proof, we proceed to the construction of the sequence $\varphi_k$. The idea is to inductively make $\varphi_{k+1}$ almost commute with the action of $G$ 
conjugated by $\tilde{\varphi}_k$ along a very small neighborhood of a large but finite part 
of the orbit of $x_0$. More precisely, let us number the points in the $G$-orbit of $x_0$ as $x_0,x_1,\ldots$
%\marginpar{agregue un comentario en el parentesis.} 
so that $x_i$ is no further to the origin than $x_j$ in the corresponding Schreier graph 
whenever $i \leq j$. (Since $x_0$ has free orbit, this Schreier graph is actually isomorphic 
to the Cayley graph of $G$.) We let $d_i$ the graph distance between $x_i$ and the origin. 
Denote $x_i^k := \tilde{\varphi}_{k-1} (x_i)$, where by definition $x_i^1=x_i$. 
Assume that the support of $\varphi_k$ consists of a collection 
of disjoint intervals $I_0^k, \ldots, I_{\ell_k}^k$ for a certain $\ell_k$, that are disjoint from $\{y_0\}$, and such that 
each $I_i^k$ is centered at $x_i^k$ and $\tilde{\varphi}_k (x_i^k) = x_i^k$ for each $i \leq \ell_k$. 
Then $\varphi_{k+1}$ will be defined so that its support consist of a collection of intervals 
$I_0^{k+1},\ldots, I_{\ell_{k+1}}^{k+1}$ disjoint from $\{y_0 \}$, where $I_j^{k+1}$ is 
centered at $x_j^{k+1}$ so that $I_j^{k+1}$ is a subset of (but much smaller than) 
$I_j^k$ for $j \leq \ell_k$, and where $\ell_{k+1}$ is to be chosen. (In particular, we will 
have $x_j^k = x_j^{k+1}$ for $j \leq \ell_k$, yet for larger $j$ these points may differ.) 
Besides, we will ask 
%\marginpar{No entiendo, para mi $x_i^k=x_i$ para todo $k$}
that $c (I_i^{k+1}) = I_j^{k+1}$ for every generator $c$ sending $x_i$ into $x_j$ 
for some $0 \leq i \leq j \leq \ell_{k+1}$. If the $I_j^{k+1}$ are chosen small enough, this 
will certainly ensure condition (i) above. The second half of condition (iv) will also be 
ensured by this property. Moreover, condition (iii) follows from the construction. 
The most involved issue concerns condition (ii).
%\marginpar{hay que elegir desde $\ell_k$ a $\ell_{k+1}$ no?} 
%Moreover, for each $x_i,x_j$ so that there is a generator $c$ 
%sending $x_i$ into $x_j$ and $i \leq j$, we assume that \hspace{0.1cm}
%$\tilde{\varphi}_k^{-1} c \tilde{\varphi}_k (I_i^{k+1}) \subset I_{j}^{k+1}$. 
%\marginpar{esto es lo mismo que $c \, I_i^{k+1} \subset I_{j}^{k+1}$}.

To deal with property (ii) (and properly define the diffeomorphism $\varphi_{k+1}$), let 
$$\bar{D} := (k+1) / \! \min D \tilde{\varphi}_k (y)$$  
and $\varepsilon = 1 / M_k 2^{k+2}$, where
\begin{equation}\label{d:epsilon}
M_k := \frac{\sup D \tilde{\varphi}_k (y)}{\inf D \tilde{\varphi}_k (y)} \cdot \sup Dc (y).
\end{equation}
Fix $\delta>0$ of the form
%\marginpar{esta frase esta muy confusa}
\hspace{0.05cm}  $\log \big( [k+1] / \min D \tilde{\varphi}_k(y) \big) /d$, \hspace{0.05cm} where $d$ is a large-enough integer, 
with $\delta$  associated in the framework of Lemma \ref{l:familias} to $\bar{D}$ and $\varepsilon$ above. 
Let $\ell_{k+1}$ be such that all points at distance $< d$ to the origin in the orbit graph appear in $x_0,\ldots,x_{\ell_{k+1}}$. 
Assume for a while that each generator $c$ of $G$ is affine on neighborhoods of $x_i^k$ for $i \leq \ell_{k+1}$, 
and further reduce the corresponding intervals $I_{i}^{k+1}$ so that they fit into these neighborhoods. Then let 
$\varphi_{k+1}$ be the diffeomorphism whose restriction to $I_i^{k+1}$ coincides with $\varphi_{I_i^{k+1}}^{D}$, 
where $D = D(i) = (k+1)^{1-d_i \delta}$, and which is the identity outside these intervals. We claim that this 
choice accomplish our needs provided the lengths of the intervals $I_{i}^{k+1}$ are very small. To see this, notice that
$$D (\tilde{\varphi}_k c \tilde{\varphi}^{-1}_k)(x) 
= 
\frac{D \tilde{\varphi}_k (c \tilde{\varphi}^{-1}_k (x))}{D \tilde{\varphi}_k (\tilde{\varphi}^{-1}_k (x))} Dc (\tilde{\varphi}^{-1}_k (x))$$
and
$$D (\tilde{\varphi}_{k+1} c \tilde{\varphi}^{-1}_{k+1})(x) 
= 
\frac{D \tilde{\varphi}_k (\varphi_{k+1} c \varphi_{k+1}^{-1} \tilde{\varphi}^{-1}_k (x))}{D \tilde{\varphi}_k (\tilde{\varphi}^{-1}_k (x))} 
\frac{D \varphi_{k+1} (c \varphi_{k+1}^{-1} \tilde{\varphi}_k^{-1} (x))}{D \varphi_{k+1} (\varphi_{k+1}^{-1} \tilde{\varphi}_k^{-1}(x))}
Dc (\varphi_{k+1}^{-1} \tilde{\varphi}^{-1}_k (x)).$$
Certainly, outside the intervals $\tilde{\varphi}_k (I_i^{k+1})$, both expressions coincide. Inside these intervals, they are essentially 
the same. More precisely, using the continuity of $D \tilde{\varphi}_k$ and $Dc$, and choosing $I_i^{k+1}$ sufficiently small, we may 
ensure that 
$$\left| \frac{D \tilde{\varphi}_k (c \tilde{\varphi}^{-1}_k (x))}{D \tilde{\varphi}_k (\tilde{\varphi}^{-1}_k (x))} Dc (\tilde{\varphi}^{-1}_k (x)) 
- \frac{D \tilde{\varphi}_k (\varphi_{k+1} c \varphi_{k+1}^{-1} \tilde{\varphi}^{-1}_k (x))}{D \tilde{\varphi}_k (\tilde{\varphi}^{-1}_k (x))} 
Dc \big( \varphi_{k+1}^{-1} \tilde{\varphi}^{-1}_k (x) \big) \right| 
\leq \frac{1}{2^{k+2}}.
$$
We are hence left to check that the factor 
$$\frac{D \varphi_{k+1} (c \varphi_{k+1}^{-1} \tilde{\varphi}_k^{-1} (x))}{D \varphi_{k+1} (\varphi_{k+1}^{-1} \tilde{\varphi}_k^{-1}(x))}$$
can be made close to 1 so that it lies between $1 - \varepsilon$ and $1+\varepsilon$. To do this, assume that $c(x_i) = x_j$ are so 
that $d_i < d_j$ (the case $d_i > d_j$ is analogous, and the case $d_i = d_j$ is impossible by the freeness of the orbit of $x_0$). 
Letting $y := \varphi_{k+1}^{-1} \tilde{\varphi}_k^{-1}(x) \in I_j^{k+1},$ we see that the expression above equals 
$$ \frac{D \varphi_{k+1} (c(y))}{D \varphi_{k+1}(y)} = \frac{D \varphi_{I_j^{k+1}}^{D(j)} (c(y))} {D \varphi_{I_i^{k+1}}^{D(i)}(y)}.$$
Now, as $c$ is affine on $I_i^{k+1}$ and sends this interval into $I_j^{k+1}$, it also transforms by conjugacy the maps 
$\varphi_{I_i^{k+1}}^D$ into $\varphi_{I_j^{k+1}}^{D}$. Hence, Lemma \ref{l:familias} yields immediately that the last 
expression lies between $1-\varepsilon$ and $1+ \varepsilon$, as desired. This yields condition (ii) of the Proposition. 
The first half of condition (iv) follows from the last step of this construction.

It remains the general case where $c$ is not affine in neighborhoods of the points $x_i^k$, with $i \leq \ell_{k+1}$. However, as the 
orbit of $x_0$ is free, we may reduce to the affine case by smooth conjugacy. Indeed, after a change of coordinates $\psi_k$, we may 
apply the arguments above, and latter come back to the original coordinates. As the distortion of this change of coordinates only 
depends on the original action and $\varphi_1, \ldots, \varphi_k$, we may keep a good control of derivatives during this process. 
More precisely, the distortion 
$$\frac{\sup D \psi_k (y)}{\inf D \psi_k (y)}$$
can be made arbitrarily close to 1 provided the lengths of the intervals $I_j^{k+1}$ are chosen small enough. 
This allows to conclude the proof.
$\hfill\square$

\vsp\vsp\vsp\vsp

Next, we deal with the $C^r$ case, where $r \geq 2$. 
%(We do not assume $r$ to be an integer.)

\vsp

\begin{prop} {\em Let $G$ be a group of the form $\Z \ltimes_A H$, where $A \in GL_d(\Q)$
has no eigenvalue of norm 1 and $rank_{\Q}(H) = d$. Then for all $r \geq 2$, every faithful
action of $G$ by $C^r$ diffeomorphisms of $[0,1]$ with no global fixed point in $(0,1)$ is
conjugate to an affine action by a homeomorphism that restricted to $(0,1)$ is a $C^r$
diffeomorphism.}
\end{prop}

\noindent{\em Proof:} 
We know from Theorem \ref{thm main} that the action is conjugate to an affine 
action via a homeomorphism $\varphi$. The image of $H$ is a subgroup of the 
group of translations which is necessarily dense; otherwise, $H$ would have 
rank 1 and $A^2$ would stabilize it pointwise, thus contradicting hyperbolicity.  
As $g$ is assumed to be $C^r$, $r\geq 2$, and has no fixed point in $(0,1)$, 
Szekeres' theorem implies that the restrictions of $g$ to $[0,1)$ and $(0,1]$
are the time-one map of the flows of vector fields $\mathcal{X_-}$ and 
$\mathcal{X_+}$, respectively, that are $C^1$ on their domains and 
$C^{r-1}$ at the interior. Futhermore, Kopell's lemma implies that the 
$C^1$ centralizer of $g$ is contained in the intersection of the flows of 
$\mathcal{X_-}$ and $\mathcal{X_+}$. Therefore, the flows coincide for 
a dense subset of times, hence $\mathcal{X_-} = \mathcal{X_+}$ on 
$(0,1)$. We denote this vector field by $\mathcal{X}$ and we call it 
the Szekeres vector field associated to $H$. 
%(and any nontrivial element $b \!\in\! H$). 
(See \cite[\S 4.1.3]{navas-book} for the details.) 

The homeomorphism $\varphi$ must send this flow into that
of the translations. Since $\mathcal{X}$ is of class $C^{r-1}$ on $(0,1)$, we have
that $\varphi$ is a $C^{r-1}$ diffeomorphism of $(0,1)$. To see that $\varphi$ 
is actually a $C^r$ diffeomorphism, we use Theorem \ref{thm multiplier}, which 
says that the interior fixed point $x_0$ of the element $a$ is hyperbolic. Indeed, this 
implies that $\varphi$ is a $C^1$ diffeomorphism that conjugates two germs of 
hyperbolic diffeomorphisms. By a well-known application of (the sharp version of) 
Sternberg's linearization theorem, such a diffeomorphism has to be of class $C^r$ 
in a neighborhood of $x_0$ (see \cite[Corollary 3.6.3]{navas-book}). Since the 
action is minimal on $(0,1)$ due to Proposition \ref{prop conjugacion}, this 
easily implies that $\varphi$ is of class $C^r$ on the whole open interval. 
$\hfill\square$

\begin{rem} Let $G$ be a group of the form $\Z \ltimes_A H$, where $A \in GL_d(\Q)$
has no eigenvalue of norm 1 and $rank_{\Q}(H) = d$.  Assume that $G$ acts by $C^r$ 
diffeomorphisms of $[0,1]$ with no global fixed point in $(0,1)$ and that $r \geq 1$. Then 
every bi-Lipschitz conjugacy of $G$ into an affine group is $C^r$ on $(0,1)$. Indeed, this 
essentially follows (with minor modifications) from \cite[\S 3.6]{navas-book}.
\end{rem}

%%%%%%%%%%%%%%%%%%%%%%%%%%%%%%%%%%%%%%%%%%%%%%%%%%%%%%%%%%%%%%%%

\section{Examples involving non-hyperbolic matrices}
\label{examples-non-hyp}

%%%%%%%%%%%%%%%%%%%%%%%%%%%%%%%%%%%%%%%%%%%%%%%%%%%%%%%%%%%%%%%%

We next consider the situation where $A \in GL_d(\Q)$ has some eigenvalues of
modulus $=1$ and some others of modulus $\neq 1$. Our goal is to prove Theorem
\ref{actuan}, according to which the group $\Z \ltimes_A \Q^d$ has an action by
$C^1$ diffeomorphisms of the closed interval that is not semiconjugate to an
affine action provided $A$ is irreducible. In particular, this is the case for the matrix
$$A := \left( \begin{array}{cccc} 0 & 0 & 0 & -1\\ 1 & 0 & 0 & -4\\ 0 & 1 & 0 & -4\\
0 & 0 & 1 & -4 \end{array} \right) \in SL_4 (\Z).$$
Indeed, $A$ has characteristic polynomial $p(x)=x^4+4x^3+4x^2+4x+1=p_1(x) p_2(x)$,
where $p_1 (x):=x^2+(2+\sqrt{2})x+1$ and $p_2 (x):=x^2+(2-\sqrt{2})x+1$. Notice that
$p(x)$ has no rational root, neither a decomposition into two polynomial of rational
coefficients of degree two; hence, it is irreducible over $\Q$. Moreover, the roots
$\lambda$ and $1/\lambda$ of $p_1$ are real numbers of modulus different
from $1$, while the roots $w$, $\overline{w}$ of $p_2$ are complex numbers
of modulus $1$, where
$$w = \frac{\sqrt{2}-2+i\sqrt{4\sqrt{2}-2}}{2},
\qquad
\lambda = \frac{-\sqrt{2}-2 + \sqrt{4\sqrt{2}+2}}{2}.$$

\vsp

Given any $A \in GL_d(\Q)$, we begin by constructing an action of
$G := \Z \ltimes_A \Q^d$ by homeomorphisms of the interval that is not 
semiconjugate to an affine action. To do this, we consider a decomposition 
$[0,1] = \overline{\bigcup_{k \in \mathbb{Z}} I_k}$, where the $I_k$'s are 
open intervals disposed on $[0,1]$ in an ordered way and such that the 
right endpoint of $I_k$ coincides with the left endpoint of $I_{k+1}$, for 
all $k \in \Z$. Let $f$ be a homeomorphism of $[0,1]$ sending each $I_{k}$
into $I_{k+1}$.  For each $(t_1,\ldots,t_d)\!\in\!\mathbb{Q}^d$ and $k \in \Z$,
denote
$$(t_{1,k},\ldots,t_{d,k}) := A^k (t_1,\ldots,t_d).$$
Let $\xi^t$ be a nontrivial topological flow on $I_0$. Next, fix
$(s_1,\ldots,s_d) \!\in\! \R^d$, and for each $(t_1,\ldots,t_d)~\in~\Q^d$,
define $g:=g_{(t_1,\ldots,t_d)}$ on $I_0$ by $g|_{I_0} := \xi^{\sum_i s_i t_i}$.
Extend $g$ to the whole interval by letting
\begin{equation}\label{extension}
g \big|_{I_{-k}} = f^{-k} \circ \xi^{\sum_i s_i t_{i,k}} \big|_{I_0} \circ f^{k}.
\end{equation}
It is not hard to see that the correspondences
$a \mapsto f,$ $(t_1,\ldots,t_d) \mapsto g_{(t_1,\ldots,t_d)},$
define a representation of $G$, where $a$ stands for the generator
of the $\Z$-factor of $G$.

\vsp
\vsp

\begin{lem} {\em If $A$ is $\Q$-irreducible and $(s_1,\ldots,s_d)$
is nonzero, then the action constructed above is faithful.}
\end{lem}

\noindent{\em Proof:} Denote by $b_1,\ldots,b_d$ the canonical basis of $H:=\Q^d$.
We need to show that for a given nontrivial $b := b_1^{t_1} \cdots b_d^{t_d} \in H$,
the associated map $g := g_{(t_1, \ldots, t_d)}$ acts nontrivially on $[0,1]$. Assume
otherwise. Then according to (\ref{extension}),  for all $k \in \mathbb{Z}$,
$$0 = \sum_i s_i t_{i,k} = \big\langle (s_1,\ldots,s_d), A^k (t_1,\ldots,t_d) \big\rangle.$$
As a consequence, the $\Q$-span of $ A^k (t_1,\ldots,t_d), k \in \mathbb{Z}$, is a
$\Q$-invariant subspace orthogonal to $(s_1,\ldots,s_d)$. However, as $A$ is
$\Q$-irreducible, the only possibility is $(t_1,\ldots,t_d) = 0$, which implies
that $b$ is the trivial element in $H$.
$\hfill\square$

\vsp
\vsp
\vsp

Assume next that $A$ is not hyperbolic. Associated to the transpose matrix
$A^T$, there is a decomposition $\R^d = E^s \oplus E^u \oplus E^c$ into stable,
unstable, and central subspaces, respectively. The space
$E^c$ necessarily contains a subspace $E^c_*$ of dimension 1 or
2 that is completely invariant under $A^T$ and such that for each nontrivial vector
therein, all vectors in its orbit under $A^T$ have the same norm. Our goal is to prove

\vsp

\begin{prop}{\em If $(s_1,\ldots,s_d)$ belongs to $E^c_*$,
then the action above is $C^1$ smoothable.}
\end{prop}

\vsp

This will follow almost directly from the next

\vsp

\begin{prop} \label{the-construction}
{\em The map $f$ and the subintervals $I_k$ of the preceding construction
can be taken so that $f$ is a
$C^1$ diffeomorphism that commutes with a $C^1$ vector field
whose
support in $(0,1)$ is nontrivial and contained in the union of the interior of the $I_k$'s.}
\end{prop}

\vsp

Using $f$ and the vector field
above, we may perform the
construction taking $\xi^t$ as being the flow associated to it. Indeed,
since the vector field
%$\mathcal{X}$
is $C^1$ on the whole interval, equation (\ref{extension}) implies
that for a given $(t_1,\ldots,t_d)$, the corresponding $g_{(t_1,\ldots,t_d)}$ is a
$C^1$ diffeomorphism provided the expressions $\sum_i s_i t_{i,k}$ remain
uniformly bounded on $k$. However, as $(s_1,\ldots,s_d)$ belongs to $E^c_*$,
this is always the case, because
$$\sum_i s_i t_{i,k}
= \big\langle (s_1,\ldots,s_d), A^k (t_1,\ldots,t_d) \big\rangle
= \big\langle (A^T)^k (s_1,\ldots,s_d), (t_1,\ldots,t_d) \big\rangle$$
and $\{ (A^T)^k (s_1,\ldots,s_d), k \in \mathbb{Z}\}$ is a bounded
subset of $\R^d$.

\vsp\vsp\vsp

To conclude the proof of Theorem \ref{actuan}, we need to show Proposition
\ref{the-construction}. Although at this point we could refer to the classical 
construction of Pixton \cite{pixton}, we prefer to give a simpler argument that 
decomposes into two elementary parts given by the next lemmas.

\vsp

\begin{lem} {\em There exists a vector field $\mathcal{X}_0$ on $[0,1]$ with compact
support in $(0,1)$ and a sequence $(\varphi_k)$ of $C^{\infty}$ diffeomorphisms of $[0,1]$
with compact support inside $(0,1)$ that converges to the identity in the $C^1$ topology
and such that the diffeomorphisms $\tilde{\varphi}_k := \varphi_k \circ \cdots \circ \varphi_1$
satisfy \hspace{0.02cm} $(\tilde{\varphi}_k)_* (\mathcal{X}_0) = t_k \mathcal{X}_0$ 
\hspace{0.02cm} for a certain sequence $(t_k)$ of positive numbers converging to zero.}
\end{lem}

\noindent{\em Proof:} Sart with the flow of translations on the real line and the
corresponding (constant) vector field. Any two positive times of this flow are smoothly
conjugate by appropriate affine transformations. Now, map the real-line into the interval
by a projective map. This yields the desired vector field and diffeomorphisms, except for
that the supports are not contained in $(0,1)$. To achieve this, just start by performing
the Muller-Tsuboi trick ({\em c.f.} Lemma \ref{MT-trick}) in order to make everything 
flat at the endpoints, then extend everything trivially in both directions by slightly 
enlarging the interval, and finally renormalize the resulting interval into $[0,1]$. 
$\hfill\square$

\vsp\vsp\vsp

Given a diffeomorphism $\varphi$ of (resp., vector field $\mathcal{X}$ on) an interval $I$,
we denote by $\varphi^{\vee}$ (resp., $\mathcal{X}^{\vee}$) the diffeomorphism of (resp., 
vector field on) $[0,1]$ obtained after conjugacy (resp., push forward) by the unique affine map
sending $I$ into $[0,1]$. Proposition \ref{the-construction} is a direct consequence of the next

\vsp

\begin{lem} {\em There exists a $C^1$ diffeomorphism $f$ of $[0,1]$ fixing
only the endpoints (with the origin as a repelling fixed point) as well as a
$C^1$ vector field $\mathcal{Y}$ on $[0,1]$ such that \hspace{0.01cm}
\noindent $f_* (\mathcal{Y}) = \mathcal{Y}$ \hspace{0.01cm} and so that for a 
certain $x_0 \in (0,1)$, we have $(\mathcal{Y}|_{[x_0,f(x_0)]})^{\vee} = \mathcal{X}_0$.}
\end{lem}

\noindent{\em Proof:}  Start with a $C^{\infty}$ diffeomorphism $g$ of $[0,1]$
that has no fixed point at the interior, and has the origin as a repelling fixed point.
Fix any $x_0 \in (0,1)$, and let $\mathcal{Z}$ be a vector field on $[x_0,g(x_0)]$
such that $\mathcal{Z}^{\vee} = \mathcal{X}_0$. A moment's reflexion shows that
this construction can be performed so that $g$ is affine close to each endpoint.
%except for a little interval contained in one of the components of
%$[x_0,f(x_0)] \setminus \supp (\mathcal{Z})$.

For each $k \in \Z$, let $I_k := g^k ( [x_0,f(x_0)] )$.
Let $\varphi_k^{\wedge}$ be a diffeomorphism of $I_k$ into itself such that
$(\varphi_k^{\wedge})^{\vee} = \varphi_k$. Now let $f$ be defined by letting 
$f \big|_{I_{|k|}} := \varphi_{|k|}^{\wedge} \circ g\big|_{I_{|k|}}$. Extend $\mathcal{Z}$ 
to the whole interval $[0,1]$ by making it commute with $g$. Finally, define $\mathcal{Y}$ 
by letting $\mathcal{Y} \big|_{I_{|k|}} := t_{|k|} \mathcal{Z} \big|_{I_{|k|}}$ for every 
$k \in \Z$. One easily checks that $f$ and $\mathcal{Y}$ satisfy the desired properties.
$\hfill\square$

\vsp\vsp\vsp

To close this section, we remark that similar ideas yield to faithful actions by
$C^1$ circle diffeomorphisms without finite orbits for the groups considered
here. Indeed, it suffices to consider $f$ as being a Denjoy counter-example
and then proceed as before along the intervals $I_k := f^{k}(I)$, where $I$
is a connected component of the complement of the exceptional minimal
set of $f$. We leave the details of this construction to the reader.

%%%%%%%%%%%%%%%%%%%%%%%%%%%%%%%%%%%%%%%%%%%%

\section{Actions on the circle}

%%%%%%%%%%%%%%%%%%%%%%%%%%%%%%%%%%%%%%%%%%%%

Recall the next folklore (and elementary) result:
For every group of circle homeomorphisms, one of the next three possibilities holds:

\vsp

\noindent (i) there is a finite orbit,

\vsp

\noindent (ii) all orbits are dense,

\vsp

\noindent (iii) there is a unique minimal invariant closed set that is homeomorphic 
to the Cantor set. (This is usually called an {\em exceptional minimal set}.)

\vsp

\noindent Moreover, a result of Margulis states that in case of a minimal action,
either the group is Abelian and conjugate to a group of rotations, or it contains
free subgroups in two generators. (See \cite[Chapter 2]{navas-book} for all of this.)

Assume next that a non-Abelian, solvable group acts faithfully by circle homeomorphisms. By 
the preceding discussion, such an action cannot be minimal. As we next show, it can admit an 
exceptional minimal set.  For concreteness, we consider the group $G := \Z \ltimes_A \Q^d$, with 
$A \in GL_d (\Q)$. Start with a Denjoy counter-example $g \in \mathrm{Homeo}_+(\mathrm{S}^1)$, 
that is, a circle homeomorphism of irrational rotation number that is not minimal. Let $\Lambda$ 
be the exceptional minimal set of $g$. Let $I$ be one of the connected components of 
$\mathrm{S}^1 \setminus \Lambda$, and for each $n \!\in\! \Z$, denote $I_n := g^n (I)$.
Consider any representation $\phi_I \!: \Q^d \rightarrow \mathrm{Homeo}(I)$. (Such an
action can be taken faithful just by integrating a topological flow up to rationally
independent times and associating the resulting maps to the generators of $\Q^d$.) Then
extend $\phi_I$ into $\phi: G \rightarrow \mathrm{Homeo}_+(\mathrm{S}^1)$ on the one
hand by letting $\phi(a) := g$, and on the other hand, for each $b \in H$,  letting
the restriction of $\phi(b)$ to $\mathrm{S}^1 \setminus \bigcup_n I_n$ being trivial,
and setting $\phi(h)|_{I_{n}} = g^{-n} \circ \phi_I (A^{-n} (h) ) \circ g^{n}$ for each
$n \in \Z$. It is easy to check that $\phi$ is faithful. Part of the content of Theorem
\ref{thm circulo} is that in case $A$ is hyperbolic, such an action cannot be by $C^1$ diffeomorphisms.
(Compare \cite{guelman-liousse}, where Cantwell-Conlon's argument is used to prove this for the
case of the Baumslag-Solitar group.)

\vsp

We next proceed to the proof of Theorem 
\ref{thm circulo}. Let again denote by $G$ a subgroup of
$\Z \ltimes_A \Q^d$ of the form $H \times_A \Z$, with $rank_{\Q}(H) = d$
and $A \in GL_d(\Z)$. Assume with no loss of generality that the canonical basis 
$\{b_1,\ldots,b_d\}$ of $\Q^d$ is contained in $H$ (see \S \ref{section-remark-base-canonica}), and
denote by $a$ the generator of the cyclic factor (induced by $A$). We start with the next

\vsp

\begin{lem}
{\em Suppose $A$ has no eigenvalue equal to 1. Then for every representation of
$G$ into $\mathrm{Homeo}_+(\mathrm{S}^1)$, the set \hspace{0.06cm} 
$\bigcap \Per (b_i)$ \hspace{0.01cm} of common periodic points of the 
$b_i$'s is nonempty and $G$-invariant.}
\end{lem}

\noindent{\em Proof:} Let $\rho_i\in \R/\Z$ be the rotation number of $b_i$. Since $H$
is Abelian and $a b_i a^{-1}= b_1^{\alpha_{1,i}} \cdots b_d^{\alpha_{d,i}}$, we have
$$\rho_i = \alpha_{1,i} \rho_1 + \cdots + \alpha_{d,i} \rho_d \quad (mod \ \Z).$$
If we denote $v := (\rho_1,\ldots,\rho_d)$, this yields $A^T v=v$ $(mod \ \Z^d)$.
Hence, $v\in (A^T-I)^{-1}(\Z^d)\subseteq \Q^d$. Therefore, all the rotation numbers
$\rho_i$ are rational, thus all the $b_i$'s have periodic points. Next, notice that for
every family of commuting circle homeomorphisms each of which has a fixed point,
there must be common fixed points. Indeed, they all necessarily fix the points in the
support of a common invariant probability measure. To show the invariance of 
\hspace{0.02cm} $\bigcap \Per (b_i)$, notice that $H$-invariance is obvious by 
commutativity. Next, let $x$ be fixed by $b_1^{k_1}, \ldots, b_d^{k_d}$. Take 
$q \in \N$ such that $q\alpha_{i,j}$ is an integer for all $i,j$. Then
$$a b_i^{q k_i} a^{-1} (x)
= b_1^{k_iq\alpha_{1,i}} \cdots b_d^{k_iq\alpha_{d,i}} (x) = x,$$
hence $b_i^{qk_i}a^{-1}(p) = a^{-1}(p)$. We thus conclude that $a^{-1}(x)$ is a
common periodic point of the $b_i$'s, as desired. $\hfill\square$

\vsp

\begin{lem}
{\em If $a$ has periodic points, then there exists a finite orbit for $G$.}
\end{lem}

\noindent{\em Proof:} If $a$ has periodic points, then every probability measure $\mu$
that is invariant by $a$ must be supported at these points. Since $G$ is solvable
(hence amenable), such a $\mu$ can be taken invariant by the whole group. The
points in the support of this measure must have a finite orbit. $\hfill\square$

\vsp
\vsp

Summarizing, for every faithful action of $G$ by circle homeomorphisms,
the nonexistence of a finite orbit implies that $a$ admits an exceptional minimal
set, say $\Lambda$. In what follows, we will show that this last possibility cannot
arise for representations into $\mathrm{Diff}^1_+(S^1)$ with non-Abelian image.

As the set $\bigcap \Per (b_i)$ is invariant under $a$, closed, and nonempty, we must
have $\Lambda \subseteq \bigcap \Per (b_i)$. Changing each $b_i$ by $b_i^k$
for some $k\in\N$, we may assume that the periodic points of the $b_i$'s are
actually fixed. (Observe that the map sending $b_i$ into $b_i^k$ and fixing
$a$ is an automorphism of $G$.) Given a point $x$ in the complement of
$\bigcap \Fix (b_i)$ (which is nonempty due to the hypothesis), 
denote by $I_{x}$ the connected component of the complement of
$\bigcap \Fix (b_i)$ containing $x$. Then there is an $H$-invariant measure $\mu_x$
supported on $I_{x}$ associated to which there is a translation vector $\tau_{x}$;
moreover, Lemma \ref{lem tauxconA} still holds in this context.

If $I$ is any connected component of the complement of $\bigcap \Fix (b_i)$, then 
there are points $z_1,\ldots z_d$ in $I$ such that $Db_i (z_i)=1$. Therefore, for 
every $\varepsilon>0$, there exists $\delta>0$ such that if $|I|<\delta$, then
$1 - \varepsilon \leq Db_i(z) \leq 1+\varepsilon$ holds for all $z \in I$ and 
all $i \in \{1,\ldots,d\}$. By decreasing $\delta$ if necessary, we may also
assume that
\begin{equation}\label{to-be-used}
1-\varepsilon \leq \frac{Da(y)}{Da(z)} \leq 1+\varepsilon \mbox{ for all } y,z
\mbox{ at distance } dist(z,y) \leq \delta.
\end{equation}

As $I_x$ is a wandering interval for $a$, we have that there exists
$k_0\in\N$ such that $|a^k(I_x)|<\delta$ and $|a^{-k}(I_x)|<\delta$
for all $k\geq k_0$. Together with (\ref{to-be-used}), this allows to
show the next analogue of Lemma~\ref{triconA} for the translation
vectors $\Delta(x) := \big( b_1(x)-x,\ldots,b_d(x)-x \big)$.

\vsp

\begin{lem}
{\em For every $\eta>0$, there exists $k_0 \in \N$ such that if we denote by
$y_k$ the left endpoint of $a^k(I)$ and we let $\varepsilon,\hat\varepsilon$
be defined by
$$\triangle(a^{-1}(x))=Da^{-1}(y_{-k}) \,A^T\triangle(x)+\epsilon(x),
\quad x\in I_{a^{-k}(x_0)}$$ and
$$\triangle(a(x))=Da(y_k) \,(A^T)^{-1}\triangle(x)+\hat\epsilon(x),
\quad x\in I_{a^{k}(x_0)},$$
then \hspace{0.02cm} $\| \epsilon(x) \| \leq \eta \big( \| \triangle(x) \| + \| \triangle (a^{-1}(x)) \| \big)$ \hspace{0.02cm} 
and \hspace{0.02cm} $\| \hat\epsilon(x) \| \leq \eta \big( \| \triangle(x) \| + \| \triangle (a^{-1}(x)) \| \big)$ 
\hspace{0.02cm} do hold for all $k \geq k_0$.}
\end{lem}

\vsp

Again, the normalized translation vectors $\vec{\tau}_{a^{-n}(x_0)}$  (resp., $\vec{\tau}_{a^{n}(x_0)}$)
accumulate at some $\vec{\tau} \in S^d$ (resp., $\vec{\tau}_*$) as $n \to \infty$. For
each $n \in \mathbb{Z}$, we let $x_n := a^{-n} (x_0)$, and we choose a sequence of
positive integers $n_k$ such that $\vec{\tau}_{x_{n_k}}\to \vec{\tau}$ and
$\vec{\tau}_{x_{-n_k}}\to \vec{\tau}_*$ as $k\to \infty$. With this notation, 
Lemma \ref{tauestri} remains true.

Finally, Lemma \ref{lem proyectivo} is easily adapted to this case:

\vsp

\begin{lem}\label{lem proy bis}{\em
For any neighborhood $V \subset S^{d-1}$ of $E^u \cap S^{d-1}_*$ in the unit 
sphere $S^{d-1}_* \subset \R^d$ (with the norm $\|\cdot\|_*$), there is $K_0 \in \N$
such that for all $k \geq K_0$ and all $x \in a^{-k}(I_{x_0})$ not fixed by $H$,
$$\frac{\triangle(x)}{\|\triangle(x)\|_*}\in V \large\implies
\frac{\triangle(a^{-1}(x))}{\| \triangle(a^{-1}(x)) \|_*}\in V.$$
Moreover, if $V$ is small enough, then there exists $\kappa > 1$ such that
$$\frac{\triangle(x)}{\| \triangle(x) \|_*} \in V
\large\implies
\| \triangle(a^{-1}x) \|_*\geq
 \kappa \hspace{0.1cm} D a^{-1}(y_{-k}) \|\triangle(x)\|_* .$$}
\end{lem}

\vsp

Now, we may conclude as in the proof of Proposition \ref{prop no hay accione por nivels}
up to a small detail. Namely, suppose
$\vec\tau_{x_0}\notin E^s$. Then $\vec\tau\in E^u$. Using Lemmas \ref{tauestri}
and \ref{lem proy bis}, we  get for $k\geq K_0$ and all $n\in \N$,
$$ \| \triangle(x_{n+k}) \|_* \geq \kappa^n Da^{-n}(y_{-k}) \| \triangle (x_k) \|_*.$$
Now, using the fact that the growth of $Da^n$ is uniformly sub-exponential,\footnote{This
is well-known and follows from the unique ergodicity of $a$ together with that the mean
of $\log(Da)$ with resepect to the unique invariant probability measure equals zero;
see \cite[Proposition I.I, Chapitre VI]{herman}.}
we get a contradiction as $n$ goes to infinity. In the case where $\vec\tau_{x_0}\in E^s$,
we have $\vec\tau_{x_0} \notin E^s$, and we may proceed as before using $a^{-1}$
instead of $a$.

\vsp

This closes the proof of the absence of an exceptional minimal set, hence of the existence of a
finite orbit for $G$.

\vsp\vsp\vsp\vsp\vsp

\noindent{\bf Acknowledgments.} We thank L.~Arenas and A. Zeghib for useful 
discussions related to \S \ref{examples-non-hyp} and to \S \ref{sec extensions}, 
respectively. We also thank \'E.~Ghys, N.~Guelman, K.~Mann, S.~Matsumoto 
and A.~Wilkinson for their interest on this work, and the anonymous referee  
for pointing out to us an error in the original version of this work as well as 
several points to improve.

All the authors were funded by 
the Center of Dynamical Systems and Related Fields~
(Anillo Project 1103, CONICYT), and would also like to thank UCN for
the hospitality during the VIII Dynamical Systems School held at San Pedro 
de Atacama (July 2013), where this work started taking its final form. 

\noindent -- 
C. Bonatti would like to thank Chicago University for its hospitality during the stay
which started his interest on this subject. 

\noindent -- 
I. Monteverde would like to thank Univ.~of Santiago for the hospitality during his 
stay in July 2013, and acknowledges the support of PEDECIBA Matem\'atica, Uruguay. 

\noindent -- 
A. Navas would like to thank Univ.~of Bourgogne for the hospitality during different 
stages of this work, and acknowledges the support of the FONDECYT Project 1120131. 

\noindent -- 
C. Rivas acknowledges the support of the CONICYT Inserci\'on Project 79130017. 

%%%%%%%%%%%%%%%%%%%%%%%%%%%%%%%%%%%%%%%%%%%%%%%%%%%%%%%%%%%%%%%%%%%%%%%%%%%%%%%%%%%%%%%%%%%%%%%%%%%%%%%%%
%%%%%%%%%%%%%%%%%%%%%%%%%%%%%%%%%%%%%%%%%%%%%%%%%%%%%%%%%%%%%%%%%%%%%%%%%%%%%%%%%%%%%%%%%%%%%%%%%%%%%%%%%

\begin{small}

%%%%%%%%%%%%%%%%%%%%%%%%%%%%%%%%%%%%%%%%%%%%%%%%%%%%%%%%%%%%%%%%%%%%%%%%%%%%%%%%%%%%%%%%%

\vspace{0.3cm}

Christian Bonatti (bonatti@u-bourgogne.fr), Univ. de Bourgogne, Dijon, France

\vsp

Ignacio Monteverde (ignacio@cmat.edu.uy), Univ. de la Rep\'ublica, Uruguay

\vsp

Andr\'es Navas (andres.navas@usach.cl), USACH, Santiago, Chile

\vsp

Crist\'obal Rivas (cristobal.rivas@usach.cl), USACH, Santiago, Chile

\end{small}

%%%%%%%%%%%%%%%%%%%%%%%%%%%%%%%%%%%%%%%%%%%%%%%%%%%%%%%%%%%%%%%%%%%%%%%%%%%%%%%%%%%%%%%%%%
%%%%%%%%%%%%%%%%%%%%%%%%%%%%%%%%%%%%%%%%%%%%%%%%%%%%%%%%%%%%%%%%%%%%%%%%%%%%%%%%%%%%%%%%%%

\end{document}